\documentstyle[12pt]{article}

\begin{document}

\author{S.V. Ludkovsky}
\title{Stochastic processes and their spectral representations
over non-archimedean fields}
\date{25.10.2007}
\maketitle

\begin{abstract}
The article is devoted to stochastic processes with values in
finite- and infinite-dimensional vector spaces over infinite fields
$\bf K$ of zero characteristics with non-trivial non-archimedean
norms. For different types of stochastic processes controlled by
measures with values in $\bf K$ and in complete topological vector
spaces over $\bf K$ stochastic integrals are investigated. Vector
valued measures and integrals in spaces over $\bf K$ are studied.
Theorems about spectral decompositions of non-archimedean stochastic
processes are proved.

\end{abstract}

\footnote{key words and phrases: stochastic processes,
non-archimedean field, zero characteristic,
random process, linear space, stochastic integral, spectral representation}

\section{Introduction}
\par Stochastic integrals and spectral representations of
stochastic processes are widely used over the fields
of real and complex numbers
\cite{dal,feller,gihsko,hentor,shir,vatach}. If consider stochastic
processes in topological groups or metric spaces it gives some
generalization, but many specific features of topological vector
spaces and results in them may naturally be missed
\cite{partasb,partasar,hentor,vatach}. At the same time
non-archimedean analysis is being fast developed in recent years
\cite{kobl,roo,sch1,vla3,escass}. It has found applications in
non-archimedean quantum mechanics and quantum field theory
\cite{vla3,aref,djord,castro,jang}. These parts of mathematical
physics heavily depend on probability theory \cite{albkar}. Then it
is very natural in the analysis on totally disconnected topological
spaces and totally disconnected topological groups
\cite{roo,luijmms03}. Stochastic processes on such groups also
permit to investigate their isometric representations in
non-archimedean spaces.
\par Remind that non-archimedean fields $\bf K$ have
non-archimedean norms, for example, for the field of $p$-adic
numbers $\bf Q_p$, where $p>1$ is a prime number
\cite{kobl,roo,wei}. Multiplicative norms in such fields $\bf K$
satisfy the strong triangle inequality: $|x+y|\le \max (|x|, |y|)$
for each $x, y\in \bf K$. \par Besides locally compact fields we
consider also non locally compact fields. For example, the algebraic
closure of $\bf Q_p$ can be supplied with the multiplicative
non-archimedean norm and its completion relative to this norm gives
the field $\bf C_p$ of complex $p$-adic numbers. The field $\bf C_p$
is algebraically closed and complete relative to its norm
\cite{kobl}. Its valuation group $\Gamma _{\bf C_p} := \{ |z|: z\in
{\bf C_p}, z\ne 0 \} $ is isomorphic with the multiplicative group
$\{ p^x: x\in {\bf Q} \} $. There exist larger fields $\bf U_p$
being extensions of $\bf Q_p$ such that $\Gamma _{\bf U_p} = \{ p^x:
x\in {\bf R} \} $. There are known extensions with the help of the
spherical completions also, if an initial field is not such
\cite{roo,sch1,diarra,escass}.

\par Stochastic processes on spaces of functions with domains of
definition in a non-archimedean linear space and with ranges in the
field of real $\bf R$ or complex numbers $\bf C$ were considered in
works \cite{byvo,evans}-\cite{evans4,khrkoz,kochubeib}. Another
types of non-archimedean stochastic processes are possible depending
on a domain of definition, a range of values of functions, values of
measures in either the real field or a non-archimedean field
\cite{kochubeia,luijmms05,ludkhr,yasuda}. Moreover, a time parameter
may be real or non-archimedean and so on, that is a lot of problems
for investigations arise.
\par Stochastic processes with values in non-archimedean spaces
appear while their studies for non-archimedean Banach spaces,
totally disconnected topological groups and manifolds
\cite{lusmfn2006}-\cite{luijmms04}. Very great importance branching
processes in graphs also have \cite{aigner,gihsko,hentor}. For
finite or infinite graphs with finite degrees of vertices there is
possible to consider their embeddings into $p$-adic graphs, which
can be embedded into locally compact fields. Considerations of such
processes reduce to processes with values in the field $\bf Q_p$ of
$p$-adic numbers. Stochastic processes on $p$-adic graphs also have
applications in analysis of flows of information, mathematical
psychology and biology \cite{khrenb}. \par More specific features
arise, when measures with values in non-archimedean fields are
considered, so this article continuous previous works of the author
in this area \cite{lusmfn2006,luijmms05,lujmsqim,ludkhr}.
\par In this article representations of
stochastic processes with values in finite- and infinite-dimensional
vector spaces over infinite fields with non-trivial non-archimedean
norms are investigated. Below different types of stochastic
processes controlled by measures with values in non-archimedean
fields of zero characteristic and stochastic integrals are studied.
Theorems about spectral decompositions of non-archimedean stochastic
processes are proved (see, for example, \S \S 20-29, 75-82, Lemmas
27 and 29, Theorems 20, 79, 81). Moreover, special features of the
non-archimedean case are elucidated. These features arise from many
differences of the classical over $\bf R$ and $\bf C$ analysis and
the non-archimedean analysis. General constructions of the paper are
illustrated in Examples 9, 31.1, 40, 74, Theorem 41, etc., where
applications to totally disconnected topological groups are
discussed as well.

\par Some necessary facts from non-archimedean probability theory
and non-archimedean analysis are recalled that to make reading
easier (see, for example, \S \S 1-6 in section 2), as well as
developed below, when it is essential. The main results of this
paper are obtained for the first time. It is necessary to note that
in this article measures and stochastic processes with values not
only in non-archimedean fields (see section 2), but also with values
in topological linear spaces which may be infinite dimensional over
non-archimedean fields are studied (see \S \S 44-74 in section 3).

\par Stochastic processes with values in $\bf Q_p^n$ have natural
interesting applications, for which a time parameter may be either
real or $p$-adic. A random trajectory in $\bf Q_p^n$ may be
continuous relative to the non-archimedean norm in $\bf Q_p$, but
its trajectory in $\bf Q^n$ relative to the usual metric induced by
the real metric may be discontinuous. This gives new approach to
spasmodic or jump or discontinuous stochastic processes with values
in $\bf Q^n$, when the latter is considered as embedded into $\bf
R^n$.

\section{Scalar spectral functions}
\par To avoid misunderstandings we first present our notations and
definitions and recall the basic facts.

\par {\bf 1. Definitions.} Let $G$ be a completely regular totally
disconnected topological space, let also $\cal R$ be its covering
ring of subsets in $G$, $\bigcup \{ A: A\in {\cal R} \} =G$. We call
the ring separating, if for each two distinct points $x, y \in G$
there exists $A\in \cal R$ such that $x\in A$, $y\notin A$. A
subfamily ${\cal S}\subset \cal R$ is called shrinking, if an
intersection of each two elements from $\cal A$ contains an element
from $\cal A$. If $\cal A$ is a shrinking family, $f: {\cal R}\to
\bf K$, where ${\bf K}=\bf R$ or $\bf K$ is the field with the
non-archimedean norm, then it is written $\lim_{A\in \cal A}
f(A)=0$, if for each $\epsilon >0$ there exists $A_0\in \cal A$ such
that $|f(A)|<\epsilon $ for each $A\in \cal A$ with $A\subset A_0$.
\par A measure $\mu : {\cal R}\to \bf K$ is a mapping with values in
the field $\bf K$ of zero characteristic with the non-archimedean
norm satisfying the following properties:
\par $(i)$ $\mu $ is additive;
\par $(ii)$ for each $A\in \cal R$ the set $\{ \mu (B): B\in {\cal R},
A\subset B \} $ is bounded; \par $(iii)$ if $\cal A$ is the
shrinking family in $\cal R$ and $\bigcap_{A\in \cal A}A =\emptyset
$, then $\lim_{A\in \cal A} \mu (A) = 0$. \par Measures on ${\sf
Bco}(G)$ are called tight measure, where ${\sf Bco}(G)$ is the ring
of clopen (simultaneously open and closed) subsets in $G$.
\par For each $A\in \cal R$ there is defined the norm:
$\| A \| _{\mu } := \sup \{ |\mu (B)|: B\subset A, B\in {\cal R} \}
$. For functions $f: G\to X$, where $X$ is a Banach space over $\bf
K$ and $\phi : G\to [0,+ \infty )$ define the norm $\| f \|_{\phi
}:= \sup \{ |f(x)| \phi (x): x\in G \} $. \par More generally for a
complete locally $\bf K$-convex space $X$ with a family of
non-archimedean semi-norms ${\cal S} = \{ u \} $ \cite{nari} define
the family of semi-norms $\| f \|_{\phi ,u}:= \sup \{ u(f(x)) \phi
(x): x\in G \} $. Recall that a subset $V$ in $X$ is called
absolutely $\bf K$-convex or a $\bf K$-disc, if $VB+VB\subseteq V$,
where $B := \{ x\in {\bf K}: |x|\le 1 \} $. Translates $x+V$ of
absolutely $\bf K$-convex sets are called $\bf K$-convex, where
$x\in X$. A topological vector space over $\bf K$ is called $\bf
K$-convex, if it has a base of $\bf K$-convex neighborhoods of zero
(see 5.202 and 5.203 \cite{nari}). A semi-norm $u$ in $X$ is called
non-archimedean, if $u(x+y)\le \max [u(x), u(y)]$ for each $x, y\in
X$. A topological vector space $X$ over $\bf K$ is locally $\bf
K$-convex if and only if its topology is generated by a family of
non-archimedean semi-norms. Therefore, a complete $\bf K$-convex
space is the projective limit of Banach spaces over $\bf K$ (see
6.204, 6.205 and 12.202 \cite{nari}).

\par Put also $N_{\mu } (x) := \inf \{ \| U
\|_{\mu }: x\in U\in {\cal R} \} $ for each $x\in G$. If a function
$f$ is a finite linear combination over the field $\bf K$ of
characteristic functions $Ch _A$ of subsets $A\subset G$ from $\cal
R$, then it is called simple. A function $f: G\to X$ is called $\mu
$-integrable, if there exists a sequence $f_1, f_2,...$ of simple
functions such that there exists $\lim_{n\to \infty } \| f-f_n \|
_{N_{\mu },u}=0$ for each $u\in \cal S$.
\par The space $L(\mu ,X)=L(G,{\cal R},\mu ,X)$ of all $\mu
$-integrable functions with values in $X$ is $\bf K$-linear. At the
same time $\int_G\sum_{j=1}^n a_jCh _{A_j}(x)\mu (dx) :=
\sum_{j=1}^n a_j\mu (A_j)$ for simple functions extends onto $L(\mu
,X)$, where $a_j\in X$, $A_j\in \cal R$ for each $j$.
\par Put ${\cal R}_{\mu } := \{ A: A\subset G, Ch _A \in L(\mu ,{\bf K})
\} $. For $A\in {\cal R}_{\mu }$ let ${\bar \mu }(A) := \int_G\chi
_A(x)\mu (dx)$.
\par For $1\le q<\infty $ denote by
\par $ \| f \|_q := [\sup_{x\in G} |f(x)|^q N_{\mu }(x)]^{1/q}$
for a simple function $f: G\to X$, when $X$ is the Banach space, or
\par $ \| f \|_{q,u} := [\sup_{x\in G} u(f(x))^q N_{\mu }(x)]^{1/q}$
for each $u\in \cal S$, when $X$ is the complete $\bf K$-convex
space. The completion of the space of all simple functions by $ \|
* \| _q$ or by $\{ \| * \|_{q,u}: u\in {\cal S} \} $
denote by $L^q(\mu ,X)$, where $L(\mu ,X)=L^1(\mu ,X)$.

\par Let $G$ be a totally disconnected completely regular space,
let also ${\sf B_c}(G)$ be a covering ring of clopen compact subsets
in $G$, suppose that $\mu : {\sf B_c}(G)\to \bf K$ is a
finitely-additive function such that its restriction $\mu |_A$ for
each $A\in {\sf B_c}(G)$ is a measure on a separating covering ring
${\cal R}(G)|_A$, where ${\sf B_c}(G)|_A = {\cal R}|_A$, ${\cal
R}|_A := \{ E\in {\cal R}: E\subseteq A \} $.
\par A measure $\eta : {\cal R}\to \bf K$ is called absolutely
continuous relative to a measure $\mu : {\cal R}\to \bf K$, if there
exists a function $f\in L(\mu ,{\bf K})$ such that $\eta (A)=\int_G
Ch _A(x) f(x)\mu (dx)$ for each $A\in \cal R$, denote it by $\eta
\preceq \mu $. If $\eta \preceq \mu $ and $\mu \preceq \eta $, then
we say that $\eta $ and $\mu $ are equivalent $\eta \sim \mu $.
\par A $\bf K$-valued measure $P$ on ${\cal R}(X)$ we call
a probability measure if $\| X \|_P =: \| P \| =1$ and $P(X)=1$ (see
\cite{ludkhr}).
\par The following statements from the non-archimedean functional
analysis proved in \cite{roo} are useful.
\par {\bf 2. Lemma} {\it Let $\mu $ be a measure on $\cal R$. There
exists a unique function $N_{\mu }: G\to [0,\infty )$ such that
\par $(1)$ $\| Ch_A \|_{N_{\mu }} = \| A \|_{\mu }$;
\par $(2)$ if $\phi : G\to [0,\infty )$ and $ \| Ch_A \|_{\phi }
\le \| A \| _{\mu }$ for each $A\in {\cal R}$, then $\phi \le N_{\mu
}$; $N_{\mu }(x) = \inf_{x\in A, A\in \cal R} \| A \|_{\mu }$ for
each $x\in X$.}
\par {\bf 3. Theorem.} {\it Let $\mu $ be a measure on $\cal R$.
Then ${\cal R}_{\mu }$ is a covering ring of $G$ and $\bar {\mu }$
is a measure on ${\cal R}_{\mu }$ that extends $\mu $.}
\par {\bf 4. Lemma.} {\it If $\mu $ is a measure on ${\cal R}$, then
$N_{\mu }=N_{\bar \mu }$ and ${\cal R}_{\mu }={\cal R}_{\bar \mu
}$.}
\par {\bf 5. Theorem.} {\it Let $\mu $  be a measure on $\cal R$,
then $N_{\mu }$ is upper semi-continuous and for every $A\in {\cal
R}_{\mu }$ and $\epsilon >0$ the set $ \{ x\in A: N_{\mu }(x)\ge
\epsilon \} $ is ${\cal R}_{\mu }$-compact.}
\par {\bf 6. Theorem.} {\it Let $\mu $ be a measure on $\cal R$,
let also $\cal S$ be a separating covering ring of $G$ which is a
sub-ring of ${\cal R}_{\mu }$ and let $\nu $ be a restriction of
$\mu $ onto $\cal S$. Then ${\cal S}_{\nu }={\cal R}_{\mu }$ and
${\bar \nu } = {\bar \mu }$.}

\par {\bf 7. Notations and definitions.}
Let $(\Omega ,{\cal A},P)$ - be a probability space, where $\Omega $
is a space of elementary events, ${\cal A}$ is a separating covering
ring of events in $\Omega $, ${\cal R}(\Omega )\subseteq {\cal
A}\subseteq {\cal R}_P(\Omega )$, $P: {\cal A}\to {\bf K}$ is a
probability, ${\bf K}$ is a non-archimedean field of zero
characteristic, $char ({\bf K})=0$, complete relative to its
multiplicative norm, ${\bf K}\supset {\bf Q_p}$, $1<p$ is a prime
number, $\bf Q_p$ is the field of $p$-adic numbers.
\par Denote by $\xi $ a random vector (a random variable for $n=1$)
with values in $\bf K^n$ or in a linear topological space $X$ over
$\bf K$ such that it has the probability distribution $P_{\xi }(A) =
P( \{ \omega \in \Omega : \xi (\omega )\in A \} )$ for each $A\in
{\cal R}(X)$, where $\xi : \Omega \to \bf X$, $\xi $ is $({\cal
A},{\cal R}(X))$-measurable, where ${\cal R}(X)$ is a separating
covering ring of $X$ such that ${\cal R}(X)\subset {\sf Bco}(X)$,
${\sf Bco}(X)$ denotes the separating covering ring of all clopen
(simultaneously closed and open) subsets in $X$. That is, $\xi
^{-1}({\cal R}(X))\subset \cal A$. If $T$ is a set and $\xi (t)$ is
a random vector for each $t\in T$, then $\xi (t)$ is called a random
function (or stochastic function). Particularly, if $T$ is a subset
in a field, then $\xi (t)$ is called a stochastic process, while
$t\in T$ is interpreted as the time parameter.
\par As usually put $M(\xi ^k) := \int_{\Omega }\xi ^k(\omega )P(d\omega )$
for a random variable $\xi $ and $k\in \bf N$ whenever it exists.

\par Random vectors $\xi $ and $\eta $ with values in $X$ are
called independent, if $P( \{ \xi \in A, \eta \in B \} )=P( \{ \xi
\in A \} ) P( \{ \eta \in B \} )$ for each $A, B\in {\cal R}(X)$.

\par {\bf 8. Definition.} Let $ \{ \Omega , {\cal R},P \} $
be a probability space with a probability measure with values in a
non-archimedean field ${\bf K}$ complete relative to its
multiplicative norm, ${\bf K} \supset \bf Q_p$. Consider a set $G$
and a ring $J$ of its subsets. Let $\xi (A)=\xi (\omega ,A)$,
$\omega \in \Omega $, be a $\bf K$ - valued random variable for each
$A\in J$ such that
\par $(M1)$ $\xi (A)\in Y$,
$\xi (\emptyset )=0$, where $Y=L^2(\Omega ,{\cal R}, P,{\bf K})$;
\par $(M2)$ $\xi (A_1\cup A_2)= \xi (A_1) +\xi (A_2)$ $mod (P)$
for each $A_1, A_2\in J$ with $A_1\cap A_2=\emptyset $;
\par $(M3)$ $ M(\xi (A_1)\xi (A_2))= \mu (A_1\cap A_2)$;
\par $(M4)$ $M(\xi (A_1)\xi (A_2))=0$ for each $A_1\cap A_2=\emptyset $,
$A_1, A_2\in J$, that is $\xi (A_1)$ and $\xi (A_2)$ are orthogonal
random variables, where $\mu (A)\in \bf K$ for each $A, A_1, A_2\in
J$.
\par The family of random variables $ \{ \xi (A): A\in J \} $
satisfying Conditions $(M1-M4)$ we shall call the elementary
orthogonal $\bf K$-valued stochastic measure.

\par {\bf 9. Example.} If $\xi (A)$ has a zero mean value $M\xi (A)=0$
for each $A\in J$, while $\xi (A_1)$ and $\xi (A_2)$ are independent
random variables for $A_1, A_2\in J$ with $A_1\cap A_2= \emptyset $,
then they are orthogonal, since $M(\xi (A_1)\xi (A_2))=(M\xi (A_1))
(M\xi (A_2))$.

\par {\bf 10. Lemma.} {\it The function $\mu $ from Definition
8 is additive.}
\par {\bf Proof.} Since $\xi (A)\in Y=L^2(\Omega ,{\cal R}, P,{\bf K})$
for each $A\in J$, then there exists $M\xi (A) = \int_{\Omega } \xi
(\omega ,A)P(d\omega )$, since $\sup_{x\in G} |\xi (\omega
,A)|^2N_P^2(\omega )\le \sup_{x\in G} |\xi (\omega ,A)|^2N_P(\omega
)$ for the probability measure $P$ having $N_p(x)\le 1$ for each
$x\in G$, that is, $L^1(\Omega ,{\cal R}, P,{\bf K})\subset
L^2(\Omega ,{\cal R}, P,{\bf K})$.
\par  Therefore, from Conditions $(M2,M4)$
for each $A_1, A_2\in J$ with the void intersection $A_1\cap A_2 =
\emptyset $ the equalities follow:
\par $M(\xi ^2 (A_1\cup A_2)) = M[(\xi (A_1)+\xi (A_2))^2]$
\\ $= M[\xi ^2(A_1) + 2\xi (A_1)\xi (A_2)+ \xi ^2(A_2)] =
M\xi ^2(A_1) + M\xi ^2(A_2)$. In view of $(M3)$ this gives
\par $\mu (A_1\cup A_2) = \mu (A_1)+ \mu (A_2)$.

\par {\bf 11. Note.} Suppose that $\mu $ has an extension to a
measure on the separating covering ring ${\cal R}(G)$, $G$ is a
totally disconnected completely regular space, where $J \subset
{\cal R}_{\mu }(G)$.

\par {\bf 12. Definitions.} Let a random function $\xi (t)$ be
with values in a complete linear locally $\bf K$-convex space $X$
over $\bf K$, $t\in T$, where $(T,\rho )$ is a metric space with a
metric $\rho $. Then $\xi (t)$ is called stochastically continuous
at a point $t_0$, if for each $\epsilon >0$ there exists $\lim_{\rho
(t,t_0)\to 0} P( \{ u(\xi (t)-\xi (t_0)) >\epsilon \} )=0$ for each
$u\in {\cal S}$. If $\xi (t)$ is stochastically continuous at each
point of a subset $E$ in $T$, then it is called stochastically
continuous on $E$.
\par If $\lim_{R\to \infty } \sup_{t\in E} P( \{ u(\xi (t))>R \} )=0$
for each $u\in \cal S$, then a random function $\xi (t)$ is called
stochastically bounded on $E$.
\par Let $L^0({\cal R}(G),X)$ denotes the class of all step (simple)
functions $f(x) = \sum_{k=1}^mc_k Ch_{A_k}(x)$, where $c_k\in X$,
$A_k\in {\cal R}(G)$ for each $k=1,...,m\in \bf N$, $A_k\cap A_j =
\emptyset $ for each $k\ne j$. Then the non-archimedean stochastic
integral by the elementary orthogonal stochastic measure $\xi (A)$
of $f\in L^0({\cal R}(G),X)$ is defined by the formula:
\par $(SI)$ $\eta (\omega ) := \int_G f(x)\xi (\omega ,dx) := \sum_{k=1}^m
c_k \xi (\omega , A_k)$.

\par {\bf 13. Lemma.} {\it Let $f, g\in L^0({\cal R}(G),{\bf K})$, where
$f(x) = \sum_{k=1}^m c_k Ch_{A_k}(x)$ and $g(x) = \sum_{k=1}^m d_k
Ch_{A_k}(x)$, then $M(\int_Gf(x)\xi (dx) \int_Gg(y)\xi (dy))=
\sum_{k=1}^m c_kd_k \mu (A_k)$ and there exists a $\bf K$-linear
embedding of $L^0({\cal R}(G),{\bf K})$ into $L^2(\mu ,{\bf K})$.}
\par {\bf Proof.} In view of Conditions $(M1,M2)$ there exists
$\int_G f(x)\xi (\omega ,dx)\in Y = L(P)$. Since $\int_Gf(x)\xi (dx)
\int_Gg(y)\xi (dy)= \sum_{k,j = 1}^m c_kd_j \xi {A_k)\xi (A_j}$,
then $M(\int_Gf(x)\xi (dx) \int_Gg(y)\xi (dy))= \sum_{k, j =1}^n
c_kd_j M(\xi (A_k)\xi (A_j)) = \sum_{k=1}^m c_kd_k \mu (A_k)$ due to
Conditions $(M3,M4)$, since $A_k\cap A_j=\emptyset $ for each $j\ne
k$. This gives the $\bf K$-linear embedding $\theta $ of $L^0({\cal
R}(G),{\bf K})$ into $L^2(\mu ,{\bf K})$ such that $\theta
(f)=\sum_{k=1}^nc_k Ch _{A_k}(x)$ and \par $(i)$ $ \| \theta (f)
\|_2=[\max_{k=1}^m |c_k|^2\sup_{x\in A_k}N_{\mu }(x)]^{1/2}=
[\max_{k=1}^m |c_k|^2 \| A_k\|_{\mu }]^{1/2}<\infty $ due to Lemma
2.

\par {\bf 14. Note.} Denote by $L^2({\cal R}(G),{\bf K})$ the completion
of $L^0({\cal R}(G),{\bf K})$ by the norm $ \| * \|_2$ induced from
$L^2(\mu ,{\bf K})$.

\par {\bf 15. Definition.}
Let $L^0({\xi },X)$ denotes the class of all step (simple) functions
$f(x) = \sum_{k=1}^mc_k \xi (A_k)$, where $c_k\in X$, $A_k\in {\cal
R}(G)$ for each $k=1,...,m\in \bf N$, $A_k\cap A_j = \emptyset $ for
each $k\ne j$. Then the non-archimedean stochastic integral by the
elementary orthogonal stochastic measure $\xi (A)$ of $f\in L^0(\xi
,X)$ is defined by the formula:
\par $(SI)$ $\eta (\omega ) := \int_G f(x)\xi (\omega ,dx) := \sum_{k=1}^m
c_k \xi (\omega , A_k)$.

\par {\bf 16. Lemma.} {\it Let $f, g\in L^0(\xi ,{\bf K})$, where
$f(x) = \sum_{k=1}^m c_k \xi (A_k)$ and $g(x) = \sum_{k=1}^m d_k \xi
(A_k)$, then $M(\int_Gf(x)\xi (dx) \int_Gg(y)\xi (dy))= \sum_{k=1}^m
c_kd_k \mu (A_k)$ and there exists a $\bf K$-linear embedding of
$L^0({\cal R}(G),{\bf K})$ into $L^2(P,{\bf K})$. }
\par {\bf Proof.}
In view of Conditions $(M1,M2)$ there exists $\int_G f(x)\xi (\omega
,dx)\in Y = L(P,{\bf K})$. But $f$ is the step function, hence \par
$(i)$ $ \| f \|_{L^2(P,{\bf K})} = [\max_{k=1}^m |c_k|^2
\sup_{\omega \in \Omega }|\xi ^2(\omega ,A_k)| N_P(x)]^{1/2}$ \\  $=
[\max_{k=1}^m |c_k|^2 \| \xi (*,A_k) \|^2_{L^2(P)}]^{1/2}<\infty $
\\ and inevitably $f\in L^2(P)$. Thus the mapping $\psi (f) :=
\sum_{k=1}^m c_k Ch_{A_k}(x)$ gives the $\bf K$-linear embedding of
$L^0(\xi ,{\bf K})$ into $L^2(P,{\bf K})$. The second statement is
verified as in Lemma 13 due to Formulas $12,15(SI)$.
\par {\bf 17. Note.} Denote by $L^2(\xi ,{\bf K})$ the completion of
$L^0(\xi ,{\bf K})$ by the norm $ \| * \|_2$ induced from
$L^2(P,{\bf K})$.

\par {\bf 18. Corollary.} {\it The mappings 12$(SI)$ and 15$(SI)$
and Conditions $(M1-M4)$ induce an isometry between $L^2({\cal
R}(G),{\bf K})$ and $L^2(\xi )$.}
\par {\bf Proof.} The valuation group
$\Gamma _{\bf K} := \{ |z|: z\in {\bf K}, z\ne 0 \} $ is contained
in $(0,\infty )$. In view of Theorem 5, Lemma 10 and Note 11 without
loss of generality for a step function $f$ we take a representation
with $A_k\in {\cal R}(G)$ such that $ \| A_k \| _{\mu } = |\mu (A_k)
|$ for each $k=1,...,m$. The family of all such step functions is
everywhere dense in $L^2({\cal R}(G),{\bf K})$.
\par Since $M(\xi
^2(A))=\mu (A)$ for each $A\in {\cal R}(G)$, then $N_{\mu }(x) =
\inf_{A\in {\cal R}(G), x\in A} \| A \|_{\mu }$, where $ \| A\|_{\mu
} = \sup \{ |\mu (B)|: B\in {\cal R}(G), B\subset A \} = \sup \{
|M(\xi ^2 (B))|: B\in {\cal R}(G), B\subset A \}$. On the other
hand, $M(\xi ^2(B)) = \int_{\Omega }\xi ^2(\omega ,B)P(d\omega )$,
$|M(\xi ^2(B))|\le \sup_{\omega \in \Omega } |\xi ^2(\omega ,B)|
N_P(\omega )$. By our supposition $\mu $ is the measure, hence
taking a shrinking family $\cal S$ in ${\cal R}(G)$ such that
$\bigcap_{A\in \cal S} A= \{ x \} $ we get
\par $N_{\mu }(x)= \inf_{A\in {\cal R}(G), x\in A}[\sup_{B\in {\cal R
}(G), B\subset A } \sup_{\omega \in \Omega } |\xi (\omega
,B)|^2N_P(\omega )]$.\\ Thus $N_{\mu }(x)=\inf_{A\in {\cal R}(G),
x\in A}[\sup_{B\in {\cal R}(G), B\subset A } \| \xi ^2(*,B)
\|_{L^2(P)}]$ and $\| A_k \|_{\mu } =  \| \xi (*,A_k) \|^2_{L^2(P)}$
for each $k=1,...,m$ due to Lemma 2 and due to the choice $ \| A_k
\| _{\mu } = |\mu (A_k) |$ above.
\par The mapping $\psi $ from \S 16 also is $\bf K$-linear from $L^0(\xi
)$ into $L^0({\cal R}(G),{\bf K})$ such that $\psi $ is the isometry
relative to $ \| *  \|_{L^2(P)}$ and $ \|
* \|_{L^2(\mu )}$ due to Formulas 13$(i)$ and 16$(i)$ and Lemma 2.
Two spaces $L^2(P)$ and $L^2(\mu )$ are complete by their
definitions, consequently, $\psi $ has the $\bf K$-linear extension
from $L^2({\cal R}(G),{\bf K})$ onto $L^2(\xi )$ which is the
isometry between $L^2({\cal R}(G),{\bf K})$ and $L^2(\xi )$.
\par {\bf 19. Definition.} If $f\in L^2({\cal R}(G),{\bf K})$,
then put by the definition:
\par $\eta =\psi (f)=\int_G f(x)\xi (dx)$.
\par The random variable $\eta $ we call the non-archimedean
stochastic integral of the function $f$ by measure $\xi $.
\par Taking a limit in $L^2(P,X)$ we denote also by $l.i.m.$.
\par {\bf 20. Theorems.} {\it 1. For a step function
$f(x)=\sum_{k=1}^n a_k Ch_{A_k} (x)$, where $a_k\in \bf K$, $A_k\in
{\cal R}(G)$, $n=n(f)\in \bf N$, the stochastic integral is given by
the formula: \par $\eta = \int f(x)\xi (dx)= \sum_{k=1}^n a_k\xi
(A_k)$.
\par 2. For each $f, g\in L^2({\cal R}(G),{\bf K})$ there is the identity:
\par $M(\int_Gf(x)\xi (dx) \int_G g(y)\xi (dy)) = \int_Gf(x)g(x)\mu
(dx)$.
\par 3. For each $f, g\in L^2({\cal R}(G),{\bf K})$ and
$\alpha , \beta \in \bf K$
the stochastic integral is $\bf K$-linear:
\par $\int_G[\alpha f(x) + \beta g(x)] \xi (dx) = \alpha \int_Gf(x)\xi (dx)
+\beta \int_G g(x)\xi (dx)$.
\par 4. For each sequence of functions $f_n\in L^2(G, {\cal R}(G),
\mu ,{\bf K})$ such that $\lim_{n\to \infty } \| f-f_n \|_{L^2(\mu
,{\bf K})}=0$ there is exists the limit: \par $\int_Gf(x)\xi (dx) =
l.i.m._{n\to \infty } \int_G f_n(x)\xi (dx)$. \par 5. There exists
an extension of $\xi $ from $\cal R$ onto ${\cal R}_{\mu }(G)$.}
\par {\bf Proof.} Statements of (1) and (3) follow from the
consideration above. To finish the proof of (2) it is sufficient to
show that $fg\in L^1(\mu ,{\bf K})$, if $f$ and $g\in L^2(\mu ,{\bf
K})$, where $\mu $ is the measure on $G$. Since $2|f(x)g(x)|\le
|f(x)|^2 +|g(x)|^2$ for each $x\in G$, then $2\sup_{x\in G}
|f(x)g(x)|N_{\mu }(x)\le \sup_{x\in G}(|f(x)|^2 +|g(x)|^2)N_{\mu
}(x)\le \| f\|^2_{L^2(\mu )}+ \| g\|^2_{L^2(\mu )}$, consequently,
$f(x)g(x)$ is $\mu $-integrable.
\par 4. From $\lim_{n\to \infty } [\sup_{x\in G}
|f(x)-f_n(x)|^2N_{\mu }(x)] =0$ and $M[\int_G(f-f_n)(x)\xi (dx)
\int_G(f-f_n)(y)\xi (dy)]= \int_G(f-f_n)^2(x)\mu (dx)$ it follows,
that $\lim_{n\to \infty }M[(\int_G(f-f_n)(x)\xi (dx))^2]=0$, that is
$l.i.m._{n\to \infty }\int_G f_n(x)\xi (dx)= \int_G f(x)\xi (dx)$
due to Corollary 18.

\par 5. Extend now the stochastic measure $\xi $ from ${\cal R}(G)$ to
${\tilde \xi }$ on ${\cal R}_{\mu }(G)$. If $A\in {\cal R}_{\mu
}(G)$, then $Ch_A\in L(G, {\cal R}(G), \mu , {\bf K})$. Since
$Ch_A\in L(G, {\cal R}(G), \mu , {\bf K})$, then $\sup_{x\in A}
N_{\mu }(x)<\infty $. Put ${\tilde \xi }(A) := \int_G Ch_A(x)\xi
(dx)=\int_A\xi (dx)$ for each $A\in {\cal R}_{\mu }(G)$,
consequently, \par $(1)$ ${\tilde \xi }$ is defined on ${\cal
R}_{\mu }(G)$. \par Therefore, ${\tilde \xi }(A)=\xi (A)$ for each
$A\in {\cal R}(G)$. For each $A, B\in {\cal R}_{\mu }(G)$ there
exist sequences of simple functions $f_n=\sum_k
a_{k,n}Ch_{A_{k,n}}$, $g_m= \sum_l b_{l,m}Ch_{B_{l,m}}$ with
$a_{k,n}, b_{l,m}\in \bf K$, $A_{k,n}, B_{l,m}\in {\cal R}(G)$ such
that $\lim_{n\to \infty } \| Ch_A - f_n\|_{L(\mu )}=0$ and
$\lim_{m\to \infty } \| Ch_B - g_m\|_{L(\mu )}=0$. Since
$M(a_{k,n}{\tilde \xi }(A_{k,n})b_{l,m}{\tilde \xi }(B_{l,m}))
=a_{k,n}b_{l,m}\mu (A_{k,n}\cap B_{l,m})$ for each $k, n, l, m$,
then \par $(2)$ $M({\tilde \xi }(A){\tilde \xi }(B)) ={\bar \mu
}(A\cap B)$ for each $A, B\in {\cal R}_{\mu }(G)$, where ${\bar \mu
}$ is the extension of the measure $\mu $ from ${\cal R}(G)$ on
${\cal R}_{\mu }(G)$. If ${\cal S}\subset {\cal R}_{\mu }(G)$ is a
shrinking family such that $\bigcap_{A\in \cal S} A=\emptyset $,
then \par $(3)$ $l.i.m._{A\in \cal S} {\tilde \xi }(A)=0$ due to
Corollary 18, since $M[{\tilde \xi }(A)]^2 = {\bar \mu }(A)$ and
$\lim_{A\in \cal S} {\bar \mu }(A)=0$ due to Theorem 3.
\par {\bf 21. Definition.} A random function of sets satisfying
conditions 20.5$(1-3)$ is called the orthogonal stochastic measure.
\par {\bf 22. Corollary.} {Let $\xi $ and $\tilde \xi $ be as in
Theorem 20.5, then $L^2(\xi ,{\bf K})=L^2({\tilde \xi },{\bf K})$.}
\par {\bf 23. Note.} If $\xi $ is an orthogonal stochastic measure
with a structure measure $\mu $ on ${\cal R}_{\mu }(G)$ and $g\in
L^2(\mu ,{\bf K})$, then put $\rho (A) := \int_GCh_A(x)g(x)\xi (dx)$
for each $A\in {\cal R}_{\mu }(G)$ and $\nu (A):=\int_Ag^2(x)\mu
(dx)$.
\par {\bf 24. Lemma.} {\it If $f\in L^2(\nu ,{\bf K})$, then $f(x)g(x)\in
L^2 (\mu ,{\bf K})$ and $\int_Gf(x)\rho (dx)=\int_G f(x)g(x)\xi
(dx)$.}
\par {\bf Proof.} In view of Theorems 20 for each
$A, B\in {\cal R}_{\mu }(G)$ there is the equality \par $M[\rho
(A)\rho (B)]= M[\int_GCh_A(x)g(x)\xi (dx) \int_G Ch_B(y)g(y)\xi
(dy)$ \par $= \int_{A\cap B} g^2(x)\mu (dx)=\nu (A\cap B)$.\\ Since
$g\in L^2(\mu ,{\bf K})$, then $\nu $ is the measure absolutely
continuous relative to $\mu $ on ${\cal R}_{\mu }(G)$. If
$f(x)=\sum_ka_kCh_{A_k}(x)$ is a simple function with $a_k\in \bf K$
and $A_k\in {\cal R}_{\mu }(G)$, then $\int_Gf(x)\rho (dx)
=\sum_ka_k\int_GCh_{A_k}(x)g(x) \xi (dx)= \sum_k a_k\rho (A_k)=
\int_Gf(x)g(x)\xi (dx)$, since \par $\sup_{x\in G}
|f(x)g(x)|^2N_{\mu }(x)\le [\max_k |a|_k^2] \sup_{x\in G}
|g(x)|^2N_{\mu }(x)<\infty .$ \\ If $f_n$ is a fundamental sequence
of simple functions in $L^2(\nu ,{\bf K})$, then
$M[(\int_G(f_n-f_m)(x)\rho (dx))^2] =
\int_G[(f_n-f_m)(x)]^2g^2(x)\mu (dx)$, hence $f_ng$ is the
fundamental sequence in $L^2(\mu ,{\bf K})$. Therefore, there exists
\par $\lim_{n\to \infty } \int_Gf_n(x)\rho (dx) = \lim_{n\to \infty }
\int_G f_n(x)g(x)\xi (dx)$, \\ consequently, $\int_Gf(x)\rho
(dx)=\int_Gf(x)g(x)\xi (dx)$.
\par {\bf 25. Lemma.} {\it If $A\in {\cal R}_{\mu }(G)$, then
$\xi (A) = \int_G[Ch_A(x)/g(x)]\rho (dx)$.}
\par {\bf Proof.} Since $\nu (\{ x: g(x)=0 \} )=0$, then $1/g(x)$
is defined $\nu $-almost everywhere on $G$, hence
$\int_G[Ch_A(x)/g^2(x)]\nu (dx)=\int_G[g^2(x)/g^2(x)]\mu (dx)=\mu
(A)$. In view of Lemma 24 $\int_G [Ch_A(x)/g(x)] \rho (dx)
=\int_G[Ch_A(x)g(x)/g(x)]\xi (dx)=\xi (A)$.

\par {\bf 26. Notation and Remark.}
 Let $T$ be a totally disconnected Hausdorff topological space with a
separating covering ring ${\cal R}(T)$ and with a non-trivial
measure $h: {\cal R}(T)\to \bf K$. Denote by $B(X,x,R) := \{ y\in X:
\rho (x,y)\le R \} $ the ball in a metric space $(X,\rho )$ with a
metric $\rho $, $0<R<\infty $. In particular, $T$ may be either a
clopen subset in $\bf K_r$ or a segment in $\bf R$, $h$ may be a
non-trivial $\bf K$-valued measure on a separating covering ring
${\cal R}(B(T,t_0,R))$ for each $t_0\in T$ and every $0<R<\infty $,
where ${\cal R}(B(T,t_0,R_1))\subset {\cal R}(B(T,t_0,R_2))$ for
each $0<R_1<R_2<\infty $ and each $t_0\in T$, ${\bf K}\supset \bf
Q_p$, ${\bf K_r}\supset \bf Q_{p'}$, $r=p'$, $r$ and $p$ are primes,
$\bf K$ and ${\bf K_r}$ are non-archimedean fields complete relative
to their multiplicative norms. \par There exists a continuous
mapping from a clopen subset in ${\bf Q_p'}$ onto $[a,b]$ in $\bf
R$, $-\infty <a<b<\infty $ (see \cite{eng}), hence a suitable
separating covering ring ${\cal R}(T)$ and $h$ exist on $[a,b]$.
\par Recall that a measure $h: {\sf B_c}({\bf K_r})\to \bf K$ is
called the Haar measure, if $h(t+B)=h(B)$ for each clopen compact
subset $B$ in $\bf K_r$ and each $t\in \bf K_r$. The Haar measure
exists due to the Monna-Springer Theorem 8.4 \cite{roo}, when $r\ne
p$ are mutually prime, $(r,p)=1$, since $B({\bf K_r},0,R)$ is
$p$-free. For example, there can be taken a non-trivial $\bf
K$-valued Haar measure $h$ on ${\sf B_c}({\bf K_r})$ such that
$h(B({\bf K_r},0,1))=1$, but generally we do not demand, that $h$ is
a Haar measure or $r\ne p$.
\par Suppose that a function $g(t,x)$ on $T\times G$ is ${\cal
R}_h\times {\cal R}_{\mu }$-measurable and $g\in L^2(T\times G,
{\cal R}_h\times {\cal R}_{\mu }, h\times \mu ,{\bf K})$, where
${\cal R}_h := {\sf B_c}_h({\bf K_r})$.

\par {\bf 27. Lemma.} {\it The stochastic integral
\par $(1)$ $\rho (t)=\int_Gg(t,x)\xi (dx)$ is defined for
each $t\in T$ for $P$-almost all $\omega \in \Omega $ and it can be
defined such that the stochastic function $\rho (t)$ would be
measurable.}
\par {\bf Proof.} If $g(t,x)=\sum_k a_k Ch_{B_k}(t) Ch_{A_k}(x)$
is a simple function with $A_k\in {\cal R}_{\mu }$ and $B_k\in {\cal
R}_h$ and $a_k\in \bf K$ for each $k=1,...,m$, $m\in \bf N$, then
$\rho (t) = \sum_k a_kCh_{B_k}(t)\xi (A_k)$ is ${\cal R}_h\times
{\cal A}$-measurable function of variables $(t,\omega )\in T\times
\Omega $ (see also Definitions 7). For each $g\in L^2(h\times \mu
,{\bf K})$ there exists a sequence of simple functions $g_n(t,x)$
such that $\lim_{n\to \infty } \sup_{t\in T, x\in G}
|g(t,x)-g_n(t,x)|^2 N_h(t) N_{\mu }(x) = 0$.
\par Let $\rho _n(t):= \int_Gg_n(t,x)\xi (dx)$, then there exists
a stochastic function ${\tilde \rho }(t)$ such that $\lim_{n\to
\infty } \sup_{t\in T} |M[({\tilde \rho }(t)-\rho
_n(t))^2]|N_h(t)=0$. There is the equality $\int_T M[({\tilde \rho
}(t)-\rho _n(t))^2]h(dt)= \int_T\int_G[g(t,x)-g_n(t,x)]^2\mu
(dx)h(dt)$, hence $\sup_{t\in T} | M[({\tilde \rho }(t)-\rho
_n(t))^2] | N_h(t) \le \sup_{t\in T, x\in G} |g(t,x)-g_n(t,x)|^2
N_h(t) N_{\mu }(x)<\infty $ and inevitably the stochastic function
${\tilde \rho }(t)$ is ${\cal R}_h\times {\cal A}$-measurable by
$(t,\omega )\in T\times \Omega $ and it exists with the unit
probability. Thus $h( \{ t\in A: M[(\rho (t) -{\tilde \rho }(t))^2]
=0 \} )= h(A)$ for each $A\in {\cal R}(T)$.
\par Finally put $\eta (t)={\tilde \rho }(t)$ if $P( \{ \rho (t)\ne
{\tilde \rho }(t) \} )=0$, while $\eta (t)=\rho (t)$ if $P( \{ \rho
(t)\ne {\tilde \rho }(t) \} )\ne 0$. Therefore, the stochastic
function $\eta $ is ${\cal R}_h\times {\cal A}$-measurable, since
$\eta $ differs from ${\cal R}_h\times {\cal A}$-measurable function
${\tilde \rho }(t)$ on a set of zero $h\times P$-measure and $\eta $
is stochastically equivalent with $\rho $.

\par {\bf 28. Remark.}
Henceforth, due to Lemma 27 we shall suppose that the stochastic
integrals 27$(1)$ are ${\cal R}_h\times {\cal A}$-measurable.

\par {\bf 29. Lemma.} {\it If $g(t,y)$ and $z(t)$ are
${\cal R}_h\times {\cal R}_{\mu }$ and ${\cal R}_h$-measurable
functions, $g\in L^2(T\times {\bf K}, {\cal R}_h\times {\cal R}_{\mu
},h\times \mu ,{\bf K})$ and $z\in L^2(T,{\cal R}_h,h,{\bf K})$,
$\xi $ is an orthogonal stochastic measure on $({\bf K}, {\cal
R}_{\mu })$, then \par $(1)$ $\int_T z(t)\int_{\bf K}g(t,y) \xi
(dy)h(dt) = \int_{\bf K}q(y)\xi (dy)$, \par where $q(y) = \int_T
z(t)g(t,y)h(dt)$.}
\par {\bf Proof.} Since $z\in L^2(h)$ and $g\in L^2(h\times \mu )$,
then \par $(2)$ $\sup_{t\in T, y\in \bf K}
|z(t)g(t,y)|^2N_h^2(t)N_{\mu }(y)\le $\par $[\sup_{t\in T}
|z(t)|^2N_h(t)] \sup_{t\in T, y\in \bf K} |g(t,y)|^2N_h(t)N_{\mu
}(y)<\infty $. \\ Consider $g\in L^2(h\times \mu )$ and a sequence
$g_n(t,y) = \sum_k a_{k,n}Ch_{B_{k,n}}(t)Ch_{A_{k,n}}(y)$ of step
functions converging to $g$ in $L^2(h\times \mu ,{\bf K})$, where
$a_{k,n}\in \bf K$, $A_{k,n}\in {\cal R}_{\mu }$, $B_{k,n}\in {\cal
R}_h$ for each $k, n$. The mean value of the square of the left side
of Equation $(1)$ is:
\par $(3)$ $M([\int_Tz(t)\int_{\bf K}g(t,y)\xi (dy)h(dt)]^2)$ \par $=
(\int_T\int_T z(t_1)z(t_2)\int_{\bf K} g(t_1,y)g(t_2,y)\mu (dy)
h(dt_1)h(dt_2)$\par $ = \int_{\bf K} [\int_T z(t)g(t,y)h(dt)]^2\mu
(dy)$.
\par Equation $(1)$ is satisfied for step functions.
In view of $(3)$ the left an the right sides of $(1)$ are continuous
relative to taking a limit by $g_n(t,y)$ in $L^2(h\times \mu ,{\bf
K})$ in the mean square sense relative to the probability $P$ as
well as in the space $L^2(P,{\bf K})$. Since the family of step
functions is dense in $L^2(h\times \mu ,{\bf K})$, then the
statement of this Lemma follows.

\par {\bf 30. Remark and Notation.} If conditions of Lemma 29
are satisfied
for each $T = B({\bf K_r},0,R)$, or $T= [-R,R]$ respectively,
$0<R<\infty $ and if there exists
\par $\int_{\bf K_r} z(t) g(t,y)h(dt) = \lim_{R\to \infty }
\int_{B({\bf K_r},0,R)} z(t)g(t,y)h(dt)$ \\ in $L^2(\mu )$, then
\par $(1)$ $\int_{\bf K_r} z(t)\int_{\bf K}g(t,y)\xi (dy) h(dt) =
\int_{\bf K} s(y)\xi (dy)$,\par where $s(y) := \int_{\bf
K_r}z(t)g(t,y)h(dt)$. This follows from Lemma 29, since the left
side of 30$(1)$ is the limit of the left side of 29$(1)$, when $R$
tends to the infinity. In the right side of 29$(1)$ it is possible
to take the limit under the sign of the stochastic integral in the
mean square sense relative to the probability $P$ as well as in the
space $L^2(P)$.

\par Describe now the generalization of the above construction onto
$\bf K$-linear spaces $X$, which may be infinite-dimensional over
$\bf K$.

\section{Vector spectral functions}

\par {\bf 31.} Let $X$ be a complete locally $\bf K$-convex space
over an infinite field $\bf K$ of zero characteristic, $char ({\bf
K})=0$, with a non-archimedean multiplicative norm relative to which
$\bf K$ is complete. Then the space $Lin (X,X) = Lin(X)$ of all $\bf
K$-linear continuous operators $F: X\to X$ is locally $\bf K$ convex
and complete. A $\bf K$-linear continuous operator $F$ is called
compact, if for each $\epsilon >0$ there exists a finite-dimensional
over $\bf K$ vector subspace $X_{\epsilon }$ such that it has a
complement $Z_{\epsilon } := X\ominus X_{\epsilon }$ in $X$ and
$u(Fx)\le \epsilon u(x)$ for each $x\in Z_{\epsilon }$ and each
semi-norm $u$ in $X$, where $Z_{\epsilon }$ is the $\bf K$-vector
subspace in $X$ such that $Z_{\epsilon }\cap X_{\epsilon }= \{ 0 \}
$, $Z_{\epsilon }\oplus X_{\epsilon } = X$. Consider the subspace
$Lc (X)$ in $Lin (X)$ of all compact operators in $X$.

\par Let $W_1: X\to X^T$ and $W_2: Lin(X) \to Lin (X)$ be linear
isomorphisms of transposition denoted simply by $W$ such that a
restriction of $W$ on each finite-dimensional subspace $\bf K^n$ in
$X$ or $Mat_n({\bf K})$ in $Lin ({\bf K})$ gives $W(F)$ a transposed
vector or matrix, where $Mat_n({\bf K})$ denotes the $\bf K$-linear
space of all $n\times n$ matrices with entries in $\bf K$. Suppose
that there exist $\bf K$-linear continuous multiplications

\par $(T1)$ $X\times X\ni \{ a, b^T \} \mapsto (a,b)\in \bf K$ and
\par $(T2)$ $X\times X\ni \{ a^T, b \} \mapsto [a,b]\in Lc(X)$, \\
where $b^T := W(b)$.

\par {\bf 31.1. Examples.} Let $X=c_0(\alpha ,{\bf K})$ be the Banach
space consisting of vectors $x= (x_j: j\in \alpha , x_j\in {\bf K})$
such that for each $\epsilon >0$ the set $ \{ j: |x_j|>\epsilon \} $
is finite with the norm $\| x \|_{c_0} :=\sup_{j\in \alpha }|x_j|$,
where $\alpha $ is a set. Due to the Zermelo Theorem (see
\cite{eng}) as $\alpha $ it can be taken an ordinal. This Banach
space $c_0(\alpha ,{\bf K})$ has the standard basis $\{ e_j: j\in
\alpha \} $, where $e_j=(0,...,0,1,0,...)$ with $1$ in the $j$-th
place and others entries zero. This basis is orthonormal in the
non-archimedean sense \cite{roo}. \par Then for each $F\in Lin (X)$
there are $F_{i,j}\in \bf K$ such that \par $(1)$ $Fe_i = \sum_{j\in
\alpha }F_{i,j}e_j$ for each $i\in \alpha $. \\ If $F\in Lc(X)$,
then for each $\epsilon >0$ the set $\beta (F) := \{ (i,j):
|F_{i,j}|> \epsilon , i, j\in \alpha \} $ is finite, where
$X_{\epsilon } = span_{\bf K} \{ e_j: \exists (i,j) \vee (j,i)\in
\beta (F) \} $, $span_{\bf K} \{ y_j: j\in \beta \} := \{
z=a_1y_{j_1}+...+a_ky_{j_k}: a_1,...,a_k\in {\bf K}, k\in {\bf N},
j_1,..., j_k\in \beta \} $ denotes the $\bf K$-linear span of
vectors.
\par If $x$ is a row-vector, then $W(x)$ is a column-vector. If $F\in Lin
(c_0(\alpha ,{\bf K}))$, then $[W(F)]_{i,j}=F_{j,i}$ for all $i,
j\in \alpha $. Taking $\epsilon _n=p^{-n}$, $n\in \bf N$, gives that
$x$ has non-zero entries only in a countable subset $\beta
(x)\subset \alpha $ and there exists $\lim_{j\in \alpha } x_j=0$. If
$a\in X$ and $b\in X$, then $(a,b)= \sum_{j\in \alpha } a_jb_j$
converges due to the non-archimedean inequality for the norm and
$\lim_{j\in \alpha } a_jb_j =0$. If $a, b\in X$, then $[a,b]=F$ with
$F_{l,j}=a_lb_j$ for all $l, j\in \alpha $, consequently, $F\in
Lc(c_0(\alpha ,{\bf K}))$.
\par If the field $\bf K$ is spherically complete, then a Banach
space over $\bf K$ is isomorphic with $c_0(\alpha ,{\bf K})$ for
some set $\alpha $ and each closed $\bf K$-linear subspace $Z$ in
$X$ is complemented (see Theorems 5.13 and 5.16 \cite{roo}). Then
certain closed $\bf K$-linear subspaces of products of Banach spaces
$c_0(\alpha ,{\bf K})$ can serve as further examples.

\par {\bf 32. Note.} Henceforth, we shall suppose that
\par $(D)$ a complete $\bf K$-convex space $X$ (see \S 31)
has an everywhere dense linear subspace $X_0$
isomorphic with $c_0(\alpha ,{\bf K})$ such that a topology $\tau
_0$ in $X_0$ inherited from the topology $\tau $ in $X$ is weaker or
equal to that of the norm topology $\tau _c$ in $c_0(\alpha ,{\bf
K})$. \par In the particular case of $\tau _0=\tau _c$ we can take
$X=c_0(\alpha ,{\bf K})$.

\par {\bf 33. Lemma.} {\it Let $X$ be a complete locally $\bf
K$-convex space satisfying Condition 32$(D)$. Then there exists a
continuous linear mapping $Tr: Lc(X)\mapsto \bf K$.}
\par {\bf Proof.} Consider an arbitrary $F\in Lc(X)$ and a semi-norm
$u$ in $X$. If a finite-dimensional over $\bf K$ subspace
$X_{\epsilon }$ is complemented in $X$, then there exists
$X_{0,\epsilon } := X_0\cap X_{\epsilon }$, $X_0\ominus X_{\epsilon
}=X_0\cap Z_{\epsilon } =: Z_{0,\epsilon }$, where $Z_{\epsilon }=
X\ominus X_{\epsilon }$ and $X_{0,\epsilon } \cap Z_{0,\epsilon } =
\{ 0 \} $, since $X_{\epsilon } \cap Z_{\epsilon } = \{ 0 \} $. Thus
$X_0= X_{0,\epsilon }\oplus Z_{0,\epsilon }$.
\par Hence there exists the continuous compact restriction of $F$ on
$X_0$. Then for each $\epsilon >0$ there exists a finite-dimensional
over $\bf K$ subspace $X_{0,\epsilon }$ in $X_0$ such that $u(Fx)\le
\epsilon u(x)$ for each $x\in Z_{0,\epsilon }$ and each semi-norm
$u$ in $X_0$. Therefore, \par $(1)$ $\lim_{i\in \alpha } \sup_{j\in
\alpha }|F_{i,j}|=0$, \\ since the family of semi-norms $ \{ u \} $
separates points in $X$, where $\{ e_j: j\in \alpha \} $ is the
basis in $X_0$ inherited from $c_0(\alpha ,{\bf K})$ and by the
Zermelo theorem we take as $\alpha $ an ordinal (see Example 31.1).
Consequently, $\sum_{j\in \alpha }F_{j,j}$ converges in $\bf K$,
since $\bf K$ is complete relative to its non-archimedean norm. Put
\par $(2)$ $Tr F|_{X_0} := \sum_{j\in \alpha}
F_{j,j}$. \\
The space $X_{0,\epsilon }$ is isomorphic with $\bf K^m$ for some
$m\in \bf N$, where the norm in $\bf K^m$ is equivalent to that of
inherited from $c_0(\alpha ,{\bf K})$. Then each basic vector $v_k$
in $X_{0,\epsilon }$ has an expansion over $\bf K$ by the basis $\{
e_j: j\in \alpha \} $, consequently, for each $\delta >0$ there
exists a finite subset $\beta $ in $\alpha $ such that $\| v_k-y_k
\|_{c_0}<\delta $ for each $k$, where $y_k\in span_{\bf K} \{ e_j:
j\in \beta \} $. \par Thus for a suitable finite subset $\beta $ in
$\alpha $ for each $x\in X_{0,\epsilon }$ there exists $y\in
span_{\bf K} \{ e_j: j\in \beta \} $ such that $u(x-y)\le \epsilon
u(x)$, since $0\le u(ax+by)\le \max (|a|u(x), |b|u(y))$ for each $x,
y\in X$ and $a, b\in \bf K$. Therefore, $Tr: Lc(X_0)\to \bf K$ is
the continuous $\bf K$-linear mapping relative to the topology $\tau
_0$ in $X_0$ provided by the family of semi-norms $\{ u \} $ and
inevitably $Tr$ has the continuous $\bf K$-linear extension on the
completion $X$ of $X_0$ relative to the locally $\bf K$-convex
topology $\tau $ in $X$.

\par {\bf 34. Corollary.} {\it Let the conditions of Lemma 33 be
satisfied and $F\in Lc(X)$. Then $\| F|_{X_0} \|_{c_0(\alpha ,{\bf
K})}<\infty $.}
\par {\bf Proof.} In view of 33$(1)$ it follows that
$\sup_{i,j\in \alpha } |F_{i,j}| <\infty $, but $\sup_{i,j\in \alpha
}|F_{i,j}| = \| F|_{X_0} \|_{c_0(\alpha ,{\bf K})} $.

\par {\bf 35. Definition.} Suppose that for each $A\in
{\cal R}(G)$ there is a random vector $\xi (A)\in X$. Let it be
satisfying the conditions:
\par $(M1)$ $\xi (A)\in Y$,
$\xi (\emptyset )=0$, where $Y=L^2(\Omega ,{\cal R}, P,X)$;
\par $(M2)$ $\xi (A_1\cup A_2)= \xi (A_1) +\xi (A_2)$ $mod (P)$
for each $A_1, A_2\in {\cal R}(G)$ with $A_1\cap A_2=\emptyset $;
\par $(M3)$ $ M[\xi (A_1),\xi (A_2)]= \mu (A_1\cap A_2)$;
\par $(M4)$ $M[\xi (A_1),\xi (A_2)]=0$ for each $A_1\cap A_2=\emptyset $,
$A_1, A_2\in {\cal R}(G)$, that is $\xi (A_1)$ and $\xi (A_2)$ are
orthogonal random variables, where $\mu (A)\in Lc(X)$ for each $A,
A_1, A_2\in {\cal R}(G)$.
\par The family of random vectors $ \{ \xi (A): A\in {\cal R}(G) \} $
satisfying Conditions $(M1-M4)$ we shall call the (elementary)
orthogonal $X$-valued stochastic measure, the compact operator $\mu
(A)$ is called the structural operator.

\par {\bf 36. Lemma.} {\it If $A_1, A_2 \in {\cal R}(G)$,
$A_1\cap A_2=\emptyset $, then $\mu (A_1\cup A_2) = \mu (A_1)+ \mu
(A_2)$.}
\par {\bf Proof.} Generalizing the proof of Lemma 10 we get the
statement of this lemma, since the product in $X$ with values in
$Lc(X)$ is continuous and $Lc(X)$ is the locally $\bf K$-convex
space having also the structure of the algebra over $\bf K$, while
$L^1(P,X)\subset L^2(P,X)$:
\par $\mu (A_1\cup A_2) = M[\xi (A_1\cup A_2), \xi (A_1\cup A_2)] =
M[\xi (A_1) +\xi (A_2), \xi (A_1) +\xi (A_2)] = M[\xi (A_1), \xi
(A_1)] + M[\xi (A_2), \xi (A_2)] = \mu (A_1) + \mu (A_2)$, since
\par $M [\xi (A_1), \xi (A_2)] =0$ and $M [\xi (A_2), \xi (A_1)] =0$
for $A_1\cap A_2 = \emptyset $.

\par {\bf 37. Note.} Generalize Definitions 1. Let $Z$ be a
locally $\bf K$-convex space with a family of semi-norms ${\cal
S}(Z)$ defining its topology. If $\cal A$ is a shrinking family, $f:
{\cal R}\to Z$, then we shall write $\lim_{A\in {\cal A}} f(A)=0$,
if for each $\epsilon >0$ and each $u\in {\cal S}(Z)$ there exists
$A_0\in \cal A$ such that $u(f(A))<\epsilon $ for each $A\in \cal A$
with $A\subset A_0$.
\par A measure $\mu : {\cal R}\to Z$ is a mapping with values in
$Z$ satisfying the following properties:
\par $(i)$ $\mu $ is additive;
\par $(ii)$ for each $A\in \cal R$ the set $\{ \mu (B): B\in {\cal R},
A\subset B \} $ is bounded; \par $(iii)$ if $\cal A$ is the
shrinking family in $\cal R$ and $\bigcap_{A\in \cal A}A =\emptyset
$, then $\lim_{A\in \cal A} \mu (A) = 0$. \par Henceforth, we
suppose that $\mu $ has an extension to a $Lc (X)$-valued measure on
the separating covering ring ${\cal R}(G)= {\cal R}$.

\par {\bf 38. Lemma.} {\it If $X$ is a complete locally
$\bf K$-convex space
satisfying Condition 32$(D)$ and $\mu $ is as in \S 37, then there
exists the trace $Tr \mu (A)$ of $\mu $ for each $A\in {\cal R}(G)$.
Moreover, $Tr \mu $ is the $\bf K$-valued measure.}
\par {\bf Proof.} In view of Lemma 33 there exists the continuous
mapping $Tr: Lc(X)\to \bf K$, hence $Tr \mu (A) \in \bf K$ for each
$A\in {\cal R}(G)$, since $\mu (A)$ is the compact operator. Then
$\| \mu (A) \|_u := \sup_{u(x)\ne 0, x\in X, } u(\mu (A)x )/ u(x)
\le \sup_{i, j\in \alpha } |[\mu (A)]_{i,j}|= \| [\mu (A)] |_{X_0}
\|_{c_0} <\infty $ due to Corollary 34, consequently,
\par $(1)$ $|Tr \mu (A)|\le \| [\mu (A)] |_{X_0} \|_{c_0}$ \\
for every $A\in {\cal R}(G)$. Therefore, if $\cal A$ is the
shrinking family in $\cal R$ and $\bigcap_{A\in \cal A}A =\emptyset
$, then $\lim_{A\in \cal A} \mu (A) = 0$, hence $\lim_{A\in \cal A}
\| \mu (A) \|_u = 0$ for each semi-norm $u$ in $X$ and inevitably
$\lim_{A\in \cal A} Tr \mu (A) = 0$ due to Inequality $(1)$ and $\mu
(A)\in Lc(X)$ and 33$(1)$.

\par {\bf 39. Definition.} Let $\sf g$ be a complete locally
$\bf K$-convex algebra with a unit $1$ and $X$ be a complete locally
$\bf K$-convex space satisfying Condition 32$(D)$ and let
simultaneously $X$ be a unital left $\sf g$-module, where $\sf g$
also satisfies Conditions 31$(T1,T2)$. This means, that there exists
a mapping ${\sf g}\times X\to X$ satisfying conditions $(1-5)$:
\par $(1)$ $b(x_1+x_2)=bx_1+bx_2$, \par $(2)$ $(b_1+b_2)x=b_1x+b_2x$,
\par $(3)$ $b_1(b_2x)=(b_1b_2)x$, \par $(4)$ $1x=x$,
\par $(5)$ there exists a
family of consistent semi-norms ${\cal S} = \{ u \} $ in $\sf g$ and
$X$ defining their Hausdorff topologies such that $u(bx)\le u(b)
u(x)$, $u(ab)\le u(a) u(b)$,
\par $(6)$ $(ax, by) = (b^Tax,y)= (x,a^Tby)\in \bf K$ and
\par $(7)$ $[ax, by] = [a,b] [x,y]\in Lc(X)$ such that $Lc(X)$
is the left $Lc ({\sf g})$-module for each $a, b, b_1, b_2\in \sf g$
and every $x, y, x_1, x_2\in X$.
\par For each $f\in L^0({\cal R},{\sf g})$ define the stochastic
integral:
\par $(SI)$ $\eta =\int_G f(x)\xi (dx):= \sum_k a_k\xi (A_k)$,
\par where $f(x) =\sum_k a_k Ch_{A_k}(x)$, $a_k\in \sf g$.
\par By $L^0(\xi ,{\sf g})=L^0(\xi )$ denote the family of all
random vectors $\eta $ of the form $(SI)$.

\par {\bf 40. Examples.} Consider either ${\sf g}=\bf K$ or
a subalgebra ${\sf g}$ in $Lin (X)$ and $X$ is a complete locally
$\bf K$-convex space, each semi-norm $v$ in $X$ induces the
consistent semi-norm $\| F \| _v := \sup_{x\in X, v(x)\ne 0}
v(Fx)/v(x)$ for each $F\in L(X)$. For simplicity of the notation we
can denote $\| F \|_v$ also by $v(F)$ and these semi-norms in $X$
and in $Lin (X)$ are consistent, since $v(Fx)\le \| F\|_v v(x)$,
where we distinguish $v(F)$ and $v(Fx)$, each multiple $bI$ of the
unit operator $I$ also belongs to $Lin (X)$ for $b\in \bf K$.

\par Take now a group $H$ with a $\bf K$-valued measure $\nu $
on ${\cal R}(H)$ or particularly $\nu $ on ${\sf B_c}(H)$ and $H$
may be a topological totally disconnected group such that
$\sup_{x\in H} N_{\nu }(x)=1$. Let $L^q_b(H,{\sf B_c}(H),\nu ,{\sf
b})$ be a completion of the family of all step functions $f: H\to
\sf b$ with supports in $A\in {\sf B_c}(H)$ on which $\nu |_A$ is
the measure relative to the family of all non-archimedean semi-norms
$ \| f\|_{q,b,v} = [\sup_{x\in H, y\in H} v[f(y^{-1}x)]^qN_{\nu
}(x)]^{1/q}<\infty $, where $1\le q<\infty $, $\sf b$ is a complete
locally convex algebra over $\bf K$ with a family of semi-norms $ \{
v \} $ in it and $\sf b$ satisfies Conditions 31$(T1,T2)$. Certainly
$v(xy)\le v(x) v(y)$ for each $x, y\in \sf b$ and every semi-norm
$v$. In particular, this space is defined for the measure $\nu $ on
${\cal R}(H)$. \par If $f_1, f_2\in L^1_b(H,\nu ,{\sf b})$, then
define the convolutions
\par $(1)$ $conv \{ f_1, f_2 \} := \{ f_1*f_2 \}(x)
:= \int_H f_1(y^{-1}x)f_2(x) \nu (dx)$ and
\par $(2)$ $conv [f_1, f_2] := [f_1*f_2](x)
:= \int_H [f_1(y^{-1}x), f_2(x)] \nu (dx)$ and
\par $(3)$ $conv (f_1,f_2) := (f_1*f_2)(x)
:= \int_H (f_1(y^{-1}x), f_2(x)) \nu (dx)$. \\
They are defined for simple functions. If they exist then
\par $(4)$ $\sup_{x\in H, z\in H} v[ \{ f_1*f_2 \} (z^{-1}x)]N_{\nu }(x)
\le $ \\ $\sup_{x\in H, y\in H, z\in H} v[f_1(y^{-1}z^{-1}x)]
v[f_2(z^{-1}x)] N_{\nu }(x) N_{\nu }(z^{-1}x)\le \| f_1 \|_{1,b,v}
\| f_2 \|_{1,b,v}<\infty $ and
\par $(5)$ $\sup_{x\in H, z\in H} v[[f_1*f_2](z^{-1}x)]N_{\nu }(x)
\le $ \\  $ \sup_{x\in H, y\in H, z\in H} C_v v[f_1(y^{-1}z^{-1}x)]
v[f_2(z^{-1}x)] N_{\nu }(x) N_{\nu }(z^{-1}x)\le C_v \| f_1
\|_{1,b,v} \| f_2 \|_{1,b,v}<\infty $ and
\par $(6)$ $\sup_{x\in H, z\in H} |(f_1*f_2)(z^{-1}x)| N_{\nu }(x)
\le $ \\   $ \sup_{x\in H, y\in H, z\in H} J_v v[f_1(y^{-1}z^{-1}x)]
v[f_2(z^{-1}x)] N_{\nu }(x) N_{\nu }(z^{-1}x)\le J_v \|
f_1 \|_{1,b,v} \| f_2 \|_{1,b,v}<\infty $, \\
since the mappings 31$(T1,T2)$ are continuous, where the semi-norm
in $Lin ({\sf b})$ induced by the semi-norm $v$ in $\sf b$ is also
denoted by $v$, $J_v$ and $C_v$ are finite semi-norms of the
mappings 31$(T1,T2)$ correspondingly relative to the semi-norm $v$
in $\sf g$. \par Therefore, the convolutions have the continuous
extensions on $L^1_b(H,\nu ,{\sf b}) =: X$ such that $conv \{
f_1,f_2 \} \in L^1_b(H,\nu ,{\sf b})$, $conv [f_1,f_2]\in
L^1_b(H,\nu ,Lc({\sf b}))$, $conv (f_1,f_2)\in L^1_b(H,\nu ,{\bf
K})$. The space $X$ is $\bf K$-linear and complete and it is the
algebra with the multiplication being the convolution $conv \{
f_1,f_2 \} $. If $1$ is not in this space adjoin it and we get the
complete locally $\bf K$-convex algebra with the unit $1(x)=1$ for
each $x\in H$.
\par This is the group algebra $X$ of $H$ over $\sf b$.
Particularly, we can take ${\sf b}=Mat_m({\bf K})$ also or more
general algebras as above. Then the transposition in $\sf b$ induces
it in $X=L^1_b(H,\nu ,{\sf b})$ such that \par $(f_1,f_2) := \int_H
(f_1(x), f_2(x))\nu (dx)\in \bf K$ and
\par $[f_1,f_2] := [f_1*f_2]\in Lin (X)$ can be considered
as the linear operator $F$ on $X$ such that $F\in Lin(X)$,
\par $(7)$ $Ff(x)= < [f_1*f_2]*f>(x)$ , where
\par $(8)$ $<g*f>(x) := \int_H g(y^{-1}x) f(x) \nu (dx)$\\ for each $g\in
L^1_b(H,{\cal R},\nu ,Lin ({\sf b}))$ and each $f\in X$ and every
$x\in H$. If $\sf b$ satisfies Condition 39$(1-7)$, then $X$ also
satisfies 39$(1-7)$. If $\sf b$ is the Banach algebra, then $X$ also
is the Banach algebra.

\par {\bf 41. Theorem.} {\it Let $Ff(x) = <[f_1*f_2]*f>(x)$, where
$f, f_1, f_2\in L^1_b(H,\nu ,{\sf b})=X$ and $H$ is the topological
group, ${\cal R}(H)\subset {\sf Bco}(H)$ as in Example 40. Then $F$
is the compact operator $F\in Lc (X)$ and the mapping $X^2\ni \{
f_1,f_2 \} \mapsto <[f_1*f_2]*> \in Lc (X)$ is continuous.}
\par {\bf Proof.} It was demonstrated in Example 40, that the convolution
$conv [f_1,f_2]$ is continuous from $X^2$ into $L^1_b(H,\nu ,Lc({\sf
b}))$, where $\nu $ is the $\bf K$-valued measure. The space
$L^1_b(H,\nu ,Lc({\sf b}))$ is the completion of the family of all
step functions $g(x)=\sum_k Ch_{A_k}(x) a_k$ relative to the family
of semi-norms $ \| g\|_{1,b,v} = [\sup_{x\in H, y\in H}
v[g(y^{-1}x)]N_{\nu }(x)]<\infty $ , where $\| Y \|_{v,\sf b} :=
\sup_{v(t)\ne 0, t\in \sf b} v(Yt)/v(t)$ is the semi-norm in $Lin
({\sf b})$ denoted also by $v(Y)$, $A_k\in {\cal R}(H)$, $a_k\in Lc
({\sf b})$. Therefore, it is sufficient to demonstrate that $F\in Lc
(X_s)$ and the mapping $X^2_s\ni \{ f_1,f_2 \} \mapsto [f_1*f_2]\in
Lc (X_s)$ is continuous, where $X_s := L^1_b(H,\nu ,{\bf K})$.
\par Theorem 7.12 \cite{roo} states, that $f\in L(H,{\cal R},\nu ,{\bf K})$
if and only if it has two properties: $(i)$ $f$ is ${\cal R}_{\nu
}$-continuous, $(ii)$ for every $\epsilon >0$ the set $\{ x:
|f(x)|N_{\nu }(x)\ge \epsilon \} $ is ${\cal R}_{\nu }$-compact,
hence contained in $ \{ x: N_{\nu }(x)\ge \delta \} $ for some
$\delta >0$. For vector valued functions see Theorem 56 below, which
is proved independently from \S \S 40, 41. Thus if $f_1, f_2\in
X_s$, then $f_1*f_2$ is ${\cal R}_{\nu }$-continuous and for each
$\epsilon >0$ and every semi-norm $v$ in $\sf b$ there exists
$\delta >0$ such that $\{ x: v(f_1*f_2(x)) N_{\nu ,v}(x)\ge \epsilon
\} \subset \{ x: N_{\nu ,v}(x)\ge \delta \} $, where $\{ x:
v(f_1*f_2(x)) N_{\nu ,v}(x)\ge \epsilon \} $ and $ \{ x: N_{\nu
,v}(x)\ge \delta \} $ are ${\cal R}_{\nu }$-compact hence ${\cal
R}$-compact sets. \par If $f\in L(H,{\cal R},\nu ,{\sf b})$, then
for each $\epsilon >0$ and each semi-norm $v$ in $\sf b$ and every
$x\in H$ there exists an open symmetric neighborhood $U_x$ of the
unit element $e$ in the topological group $H$ such that
$v(f(y^{-1}x)-f(x))<\epsilon $ for each $y\in U_x^3$. From the
covering $\{ xU_x: x\in  H, N_{\nu ,v}(x)\ge \delta \} $ of $ \{ x:
N_{\nu ,v}(x)\ge \delta \} $ extract a finite covering $\{
x_jU_{x_j}: j=1,...,q \} $ and take $U=\bigcap_{j=1}^q U_{x_j}$,
since ${\cal R}\subset {\sf Bco}(H)$. Then $U$ is open symmetric
$U^{-1}=U$ and $e\in U$. \par If $y\in U$ and $N_{\nu ,v}(x)\ge
\delta $, then there exists $j$ such that $x\in x_jU_{x_j}$, hence
$v(f(y^{-1}x)-f(x))\le \max (v(f(y^{-1}x)-f(y^{-1}x_j)),
v(f(y^{-1}x_j)-f(x_j)), v(f(x_j)-f(x)))<\epsilon $, since
$(y^{-1}x)(y^{-1}x_j)^{-1}\in U^3\subset U_{x_j}^3$. Consequently,
$v(f(y^{-1}z)- f(z))\le \max (v(f(y^{-1}z)-f(y^{-1}t)),
v(f(y^{-1}t)- f(t)))< \epsilon $ for each $z\in [U \{ x\in H: N_{\nu
,v}(x)\ge \delta \} ]$, where $t\in \{ x\in H: N_{\nu ,v}(x)\ge
\delta \} $ is such that $zt^{-1}\in U$, since
$(y^{-1}z)(y^{-1}t)^{-1}\in U^3$. At the same time $\{ x\in H:
N_{\nu ,v}(x)\ge \delta \}\subset [U \{ x\in H: N_{\nu ,v}(x)\ge
\delta \} ]$ and $v(f(z)) N_{\nu ,v}(z)<\epsilon $ for each $z\in
H\setminus \{ x\in H: N_{\nu ,v}(x)\ge \delta \}$. If $z\in
H\setminus [U \{ x\in H: N_{\nu ,v}(x)\ge \delta \} ]$, then $Uz\cap
\{ x\in H: N_{\nu ,v}(x)\ge \delta \} =\emptyset $. Thus for each
$f\in L(H,{\cal R},\nu ,{\sf b})$ and every $\epsilon
>0$ and each semi-norm $v$ in $\sf b$ there exists an open symmetric
neighborhood $U$ of $e$ in $H$ such that $v(f(y^{-1}x)-f(x))
<\epsilon $ for each $x\in H$ and $y\in U$.
\par Suppose that $f, f_1, f_2\in X_s$, then for each $\epsilon >0$
there exists an open symmetric neighborhood $U$ of $e$ in $H$ such
that \par $(iii)$ $v( < [ f_1*f_2 ] *f > (y^{-1}x) - < [ f_1*f_2
] *f > (x)) N_{\nu ,v}(x)<\epsilon \sup_{x\in H} v(f(x))  N_{\nu ,v}(x)$\\
for each $x\in H$ and every $y\in U$. Take the sets $\{ x: v( [
f_1*f_2 ] (x)) N_{\nu ,v}(x)\ge \epsilon \} \subset \{ x: N_{\nu
,v}(x)\ge \delta \} =: A$ for $[ f_1*f_2 ]$ as above. Decompose
$F=F_1+F_2$, where $F_1f := < [ f_1*f_2 ] *(Ch_A f) > $ and $F_2f :=
< [ f_1*f_2 ] *((1-Ch_A) f > $. Then $\| F_2\|_v \le \epsilon $.
\par In view of Inequality $(iii)$ for each $\epsilon >0$ there
exists a finite covering family of subsets $A_1,...,A_m\in {\cal R}$
in $A$, $\bigcup_{j=1}^mA_j=A$, such that for each $f\in X_s$ there
exists a simple function $g(x)=\sum_{j=1}^m b_kCh_{A_k}(x)$ for
which $\| F_1(f-g) \|_v <\epsilon $. But the family of such simple
functions $g$ is finite-dimensional over $\bf K$, consequently, $F$
is the compact operator.
\par In view of 40$(2,5)$ the mapping $X^2\ni
\{ f_1, f_2 \} \mapsto <[f_1*f_2]*> \in Lc (X)$ is continuous.

\par {\bf 42. Lemma.} {\it Let for a natural number $k$ there exists
$M\xi ^k$, where $\xi $ is a random vector with values in a locally
$\bf K$-convex algebra $\sf g$, then for each $1\le l<k$ there
exists $M\xi ^l$.}
\par {\bf Proof.} If $\xi $ is a random vector with values
in $\sf g$, then $\xi $ by the definition it is $({\cal A},{\cal
R}({\sf g}))$-measurable and for each semi-norm $u$ in $\sf g$ there
exists $\sup_{\omega \in \Omega } [u(\xi )^kN_P(\omega )]^{1/k} = \|
\xi \|_{k,u}<\infty $, where ${\cal R}({\sf g})$ is a separating
covering ring of $\sf g$ such that ${\cal R}({\sf g})\subset {\sf
Bco}({\sf g})$. Therefore, $\sup_{\omega \in \Omega } [u(\xi
)^lN_P(\omega )]^{1/l}\le \sup_{\omega \in \Omega } [u(\xi
)^kN_P(\omega )]^{1/k}$ for each $1\le l<k$, since $N_P(\omega )\le
1$ for each $\omega \in \Omega $.

\par {\bf 43. Remark.} Denote by $L^2(\xi ,{\sf g})$ the completion of
$L^0(\xi ,{\sf g})$ by the family of semi-norms \par $(1)$ $ \| f
\|_{2,P,u} := [\sup_{x\in G, \omega \in \Omega } u^2(f(x)\xi (\omega
,x)) N_P(\omega )]^{1/2}$ \\ induced from $L^2(\Omega ,{\cal
A},P,X)$, where $L^0(\xi ,{\sf g})$ is the space of all step
functions $\sum_{k=1}^l a_k\xi (A_k)$, $a_k\in \sf g$, $A_k\in {\cal
R}(G)$, $A_j\cap A_k= \emptyset $ for each $k\ne j$, $l\in \bf N$.
\par {\bf 44. Definition.} Suppose that $Z$ is a locally
$\bf K$-convex space and $\mu : {\cal R}\to Z$ is a measure, where
${\cal R}={\cal R}(G)$ is a separating covering ring of a set $G$.
For each semi-norm $u$ in $Z$ and any $A\in {\cal R}(G)$ define
\par $(1)$ $\| A \|_{\mu ,u} := \sup \{ u(\mu (B)): B\in {\cal R}, B\subset
A \} $.

\par {\bf 45. Lemma.} {\it If $\mu : {\cal R}\to Z$ is a measure,
then for each semi-norm $u$ in $Z$ there exists a unique function
$N_{\mu ,u}: G\to [0,\infty )$ such that \par $(1)$ $\| Ch_A
\|_{N_{\mu ,u}} = \| A \|_{\mu ,u}$ for each $A\in {\cal R}$ (see
Definition 44);
\par $(2)$ if $\phi : G\to [0,\infty )$ and $\| Ch_A \|_{\phi }\le
\| A \|_{\mu , u}$ for all $A\in {\cal R}$, then $\phi \le N_{\mu
,u}$; moreover, \par $(3)$ $N_{\mu ,u} (x)= \inf_{A: x\in A\in {\cal
R}} \| A \|_{\mu ,u}$ for each $x\in G$.}
\par {\bf Proof.} If $N_{\mu ,u}(x)$ is defined by the formula $(3)$,
then $(2)$ is evident, since $\| Ch_A \|_{\phi }\le \| A \|_{\mu ,
u}$ for all $A\in {\cal R}$. Take $A\in {\cal R}$ and consider the
family ${\cal E} := \{ B\in {\cal R}: B\subset A, \| A\setminus B
\|_{\mu ,u}\le \| Ch_A \|_{N_{\mu ,u}} +\epsilon \} $. \par On the
other hand, from 44$(1)$ it follows that $ \| A_1\cup A_2 \|_{\mu
,u} \le \max( \| A_1 \|_{\mu ,u}, \| A_2 \|_{\mu ,u} \} $ for each
$A_1, A_2\in {\cal R}$, since $\mu $ is additive and for each
$B\subset A_1\cup A_2$ we have $B=B_1\cup B_2$, while $u(\mu
(A_1\cup A_2))\le \max \{ u(\mu (A_1\setminus A_2)), u(\mu
(A_2\setminus A_1)), u(\mu (A_1\cap A_2)) \} $, where $B_1 :=
A_1\cap B$ and $B_2 := A_2\cap B$. Therefore, $\cal E$ is the
shrinking family.
\par Then for each $x\in A$ there exists $B\in\cal R$ such that
$x\in B$ and $ \| B \|_{\mu ,u} \le N_{\mu ,u}(x)+\epsilon \le \|
Ch_B \|_{N_{\mu ,u}} +\epsilon $, consequently, $A\setminus B \in
\cal E$. Then $\bigcap_{A\in \cal E} A= \emptyset $, consequently,
there exists $B\in \cal E$ such that $ \| B\|_{\mu ,u}\le \epsilon $
and inevitably $\| A \|_{\mu ,u} \le \max \{ \| B \|_{\mu ,u}, \|
A\setminus B \|_{\mu ,u} +\epsilon \} $.

\par {\bf 46. Definitions.} Let $\mu $ be a $X$-valued measure as in
\S \S 37, 39. For $f: G\to \sf g$ consider the family of semi-norms
\par $\| f\|_{q,\mu ,u} := [\sup_{x\in G} u^q(f) N_{\mu
,u}(x)]^{1/q}$\\ whenever it exists, where $1\le q<\infty $ and $u$
is a semi-norm in $X$. Define the space $L^q(G,{\cal R},\mu ,X;{\sf
g})$ as the completion of the family of all step (simple) functions
$f$ relative to the family of semi-norms $\{ \|* \|_{q, \mu ,u}:
u\in {\cal S} \} $. If $f\in L^1(G,{\cal R},\mu ,X;{\sf g})$, then
it is called $\mu $-integrable. \par In the particular case of ${\sf
g}=\bf K$ we can omit it from the notation. When $G$, $\cal R$, $\mu
$, $X$ and $\sf g$ are specified it also can be written shortly
$L^q(\mu )$, for $q=1$ it can be omitted writing $L(\mu )$. \par A
function $f$ is called $\mu $-integrable, if $f\in L(\mu )$. Put
${\cal R}_{\mu } := \{ A\subset G: Ch_A\in L(\mu ) \} $ for the
$X$-valued measure $\mu $ and extend it by ${\bar \mu }(A) := \int_G
Ch_A(x)\mu (dx)$, since $1\in \sf g$.

\par {\bf 47. Lemma.} {\it If $\mu $ and ${\cal R}_{\mu }$ are as in
\S 46, then $A\in {\cal R}_{\mu }$ if and only if for each $\epsilon
>0$ and each semi-norm $u$ in $X$ there exists $B\in {\cal R}$ such
that $N_{\mu ,u}(x)\le \epsilon $ for every $x\in A\bigtriangleup
B$.}
\par {\bf Proof.} If for each $\epsilon >0$ and a semi-norm $u$
in $X$ there exists $B\in \cal R$ such that $\| Ch_A - Ch_B
\|_{N_{\mu ,u}}<\epsilon $, then taking a sequence $\epsilon
_n=p^{-n}$ we get, that $Ch_A\in L(\mu )$, hence $A\in {\cal R}_{\mu
}$. \par On the other hand, if $A\in {\cal R}_{\mu }$, then for each
$1>\epsilon >0$ there exists a simple function $f\in L(\mu )$ such
that $ \| Ch_A - f \|_{N_{\mu ,u}} < \epsilon $. Consider $B := \{
x: u(f(x)- 1) <1 \} $ which is in $\cal R$, $B\in \cal R$. Then $
u(f(x) -Ch_B(x)) \le \min [u(f(x)), u(f(x) - 1)] \le u(f(x) -
Ch_A(x))$ for each $x\in G$ and inevitably $\| Ch_A - Ch_B
\|_{N_{\mu ,u}}\le \max [ \| f-Ch_A \|_{N_{\mu ,u}}, \| f - Ch_B
\|_{N_{\mu ,u}}] = \| f-Ch_A \|_{N_{\mu ,u}}\le \epsilon $.

\par {\bf 48. Lemma.} {\it If $\mu : {\cal R}\to X$ satisfies Conditions
37$(i,ii)$, then 37$(iii)$ is equivalent to: \par $(iii')$ if ${\cal
A}\subset {\cal R}$ is a shrinking family and $\bigcap_{A\in \cal A}
A = \emptyset $, then $\lim_{A\in \cal A} \| A \| _{\mu ,u} =0$ for
each semi-norm $u\in \cal S$ in $X$.}
\par {\bf Proof.} Since $\| A\|_{\mu ,u} \ge u(\mu (A))$ for each
$A\in \cal R$, then $(iii')$ implies 37$(iii)$.
\par Prove now the converse statement supposing that $\mu $ satisfies
37$(i-iii)$. Suppose that $\cal A$ is a shrinking family with
$\bigcap_{A\in \cal A} A=\emptyset $. For each $\epsilon >0$ and
every $u\in \cal S$ there exist $E\in \cal A$ such that $u(\mu
(A))<\epsilon $ for each $A\in \cal A$ such that $A\subset E$. Take
a semi-norm $u\in \cal S$ in $X$. For every $A\in \cal A$ choose
$V_A\in \cal R$ such that $V_A\subset A$ and $u(\mu (V_A))>\min
(\epsilon , \| A\|_{\mu ,u} /2 )$. If $A\in \cal A$, then the family
${\cal C}:= \{ V_A\cap B: B\in {\cal A}, B\subset A \} $ is
shrinking and $\bigcap_{U\in \cal C} U= \emptyset $, hence for each
$A\in \cal A$ there exists $W_A\in \cal A$ such that $W_A\subset A$
and $u(\mu (V_A\cap W_A))< \epsilon $.
\par Take the family ${\cal V} := \{ V_A\cup W_A: A\in {\cal A},
A\subset E \} $. If $A, B\in \cal A$, then there exists $C\in \cal
A$ with $C\subset W_A\cap W_B$. Then $V_C\cup W_C\subset C\subset
W_A\subset V_A\cup W_A$, also $V_C\cup W_C\subset V_B\cup W_B$,
consequently, the family $\cal V$ is shrinking. Moreover,
$\bigcap_{C\in \cal V} C = \emptyset $. Thus there exists $A\in \cal
A$ with $A\subset E$ and $u(\mu (V_A\cup W_A))<\epsilon $, also
$u(\mu (V_A\cup W_A))<\epsilon $ and $u(\mu (V_A\cap W_A))<\epsilon
$ due to the definition of $W_A$, as well as $u(\mu (W_A))<\epsilon
$, since $W_A\in \cal A$ and $W_A\subset A\subset E$. Hence $u(\mu
(V_A)) = u(\mu (V_A\cup W_A) + \mu (V_A\cap W_A) - \mu (W_A))
<\epsilon $. Therefore, $\| A \|_{\mu ,u}\le 2\epsilon $, since
$u(\mu (V_A)) >\min (\epsilon , \| A \|_{\mu ,u}/2)$.

\par {\bf 49. Theorem.} {\it Let $\mu $ be a $X$-valued measure
on $\cal R$. Then ${\cal R}_{\mu }$ is a covering ring of $G$ and
$\bar \mu $ is a $X$-valued measure extending $\mu $.}
\par {\bf Proof.} In view of Lemma 47 ${\cal R}_{\mu }$ is a
covering ring of $G$ and $\bar \mu $ is additive. If $A\in {\cal
R}_{\mu }$, then for each $B\subset A$ such that $B\in {\cal R}$ and
each semi-norm $u$ in $X$ the inequalities are satisfied:
\par $u({\bar \mu }(B)) \le \| Ch_B \|_{N_{\mu ,u}}\le
\| Ch_A \|_{N_{\mu ,u}}<\infty $, hence $\bar \mu $ has property
37$(ii)$.
\par Consider now any shrinking family ${\cal A}\subset {\cal R}_{\mu }$
having empty intersection. For $\epsilon >0$ and a semi-norm $u$ in
$X$ take ${\cal E} := \{ B\in {\cal R}: \exists A\in {\cal A}$
$\mbox{such that}$ $A\cap G_{\epsilon ,u}= B\cap G_{\epsilon ,u} \}
$, where $G_{\epsilon ,u} := \{ x\in G: N_{\mu ,u}(x)\ge \epsilon \}
$.  Then ${\cal E}$ is shrinking. If $x\notin G_{\epsilon ,u}$, then
there exists $V\in \cal R$ such that $x\in V$ and $\epsilon >
 \| V \|_{\mu ,u} $, hence $B\setminus V \in \cal E$ for each
$B\in \cal E$, consequently, $V\cap (\bigcap_{B\in \cal E} B)=
\emptyset $ and inevitably $\bigcap_{B\in \cal E} B \subset
G_{\epsilon ,u}$. Then by the construction of $\cal E$ we get:
$\bigcap_{B\in \cal E} B = \bigcap_{B\in \cal E} B\cap G_{\epsilon
,u} = \bigcap_{A\in \cal A} A\cap G_{\epsilon ,u} = \emptyset $.
\par In view of Lemma 48 there exists $B\in \cal E$ such that
$\| B \|_{\mu ,u}<\epsilon $, hence $B\cap G_{\epsilon ,u}=\emptyset
$. Then there exists $A\in \cal A$ such that $A\cap G_{\epsilon ,u}
= B\cap G_{\epsilon ,u}=\emptyset $, consequently, $ \| A \|_{\mu
,u}<\epsilon $. Again by Lemma 48 $\lim_{A\in \cal A} {\bar \mu
}(A)=0$. Thus $\bar \mu $ is the $X$-valued measure.

\par {\bf 50. Lemma.} {\it If $\mu $ is a $X$-valued measure on
$\cal R$, then $N_{\mu ,u} = N_{{\bar \mu },u}$ for each semi-norm
$u$ in $X$. Therefore, $\| * \|_{N_{\mu ,u}} = \| * \|_{N_{{\bar \mu
},u}}$, $L({\bar \mu })=L(\mu )$, $\int_G fd{\bar \mu } =
\int_Gfd\mu $, ${\cal R}_{\bar \mu } = {\cal R}_{\mu }$.}
\par {\bf Proof.} Take a semi-norm $u$ in $X$ and a point $x\in G$
and a number $b>N_{\mu ,u}(x)$. Then there exists $A\in {\cal
R}\subset {\cal R}_{\mu }$ such that $x\in A$ and $ \| A\|_{N_{{\bar
\mu },u}} \le b$. Then for every $B\in \cal R$ such that $B\subset
A$ there are inequalities $u ({\bar \mu }(B)) \le \| B \|_{N_{{\bar
\mu },u}} \le \| A \|_{N_{{\bar \mu },u}} \le b$, hence \par $(1)$ $
\| A \| _{N_{{\bar \mu },u}} \le b$, consequently, $N_{{\bar \mu
},u}(x)\ge \inf \{ \| A\|_{{\bar \mu },u}: A\in {\cal R}_{\mu },
x\in A \} $.
\par Take now $0<d< N_{\mu ,u}(x)$ and $A\in {\cal R}_{\mu }$ with
$x\in A$. In view of Lemma 47 there exists $B\in \cal R$ such that
$N_{{\bar \mu },u}(y)\le d$ for each $y\in A\bigtriangleup B$.
Therefore, $\| B \|_{N_{{\bar \mu },u}}\ge N_{{\bar \mu },u}>d$, so
$u(\mu (E))>d$ for some $E\in \cal R$ such that $E\subset B$. As
$u(\mu (E) - {\bar \mu }(E\cap A)) =u({\bar \mu } (E\setminus A)\le
\| E\setminus A \|_{N_{{\bar \mu },u}}\le \| E\setminus B
\|_{N_{{\bar \mu },u}} \le d < u(\mu (E))$. Thus $u({\bar \mu
}(E\cap A)) = u(\mu (E))$, consequently, \par $(2)$ $\| A \|_{{\bar
\mu },u} \ge u({\bar \mu } (E\cap A)) = u(\mu (E)) >d$. Finally,
from these two Inequalities $(1,2)$ the equality $N_{{\bar \mu
},u}(x) = N_{\mu ,u}(x)$ for each $x\in G$ follows.

\par {\bf 51. Theorem.} {\it $(1)$. If $\mu $ is a measure on $\cal
R$, then $N_{\mu ,u}$ is upper semi-continuous for each semi-norm
$u\in \cal S$ in $X$, hence it is ${\cal R}_{\mu }$-upper
semi-continuous and for every $A\in {\cal R}_{\mu }$ and $\epsilon
>0$ the set $\{ x\in A: N_{\mu ,u} (x)\ge \epsilon \} $ is ${\cal
R}_{\mu }$-compact, hence it is ${\cal R}$-compact.
\par $(2)$. Conversely, let $\mu : {\cal R}\to X$ satisfies
37$(i)$ and let for every $u\in \cal S$ there exist an $\cal
R$-semi-continuous function $\phi _u: G\to [0,\infty )$ such that
$u(\mu (A))\le \sup_{x\in A} \phi _u(x)$ for each $A\in \cal R$ and
let the set $ \{ x\in A: \phi _u(x)\ge \epsilon \} $ is $\cal
R$-compact for each $\epsilon >0$. Then $\mu $ is a $X$-valued
measure and $N_{\mu ,u}(x)\le \phi _u(x)$ for each $x\in G$ and each
$u\in \cal S$.}
\par {\bf Proof.} $(1)$. Put $G_{\epsilon ,u} := \{ x\in G: N_{\mu
,u}(x)\ge \epsilon \} $, where $\epsilon >0$. Then for each $x\in
G\setminus G_{\epsilon ,u}$ there exists $A\in \cal R$ such that
$x\in A$ and $\| A \|_{\mu ,u}<\epsilon $, hence $A\subset
G\setminus G_{\epsilon ,u}$ and inevitably $G_{\epsilon ,u}$ is
$\cal R$-closed and $N_{\mu ,u}$ is $\cal R$-upper semi-continuous.
Take now $A\in {\cal R}_{\mu }$ and a covering $\cal V$ of $A\cap
G_{\epsilon ,u}$ by elements of ${\cal R}_{\mu }$. Then the sets
$A\setminus (V_1\cup ... \cup V_n\cup V)$, where $n\in \bf N$,
$V_1,...,V_n\in \cal V$, $V\in {\cal R}_{\mu }$ and $V\subset G
\setminus G_{\epsilon ,u}$, form a shrinking subfamily $\cal A$ in
${\cal R}_{\mu }$ with the empty intersection $\bigcap_{E\in \cal A}
E=\emptyset $. In accordance with Property 48$(iii')$ of an
$X$-valued measure there exist $V_1,...,V_n\in \cal V$ and $V\subset
G \setminus G_{\epsilon ,u}$ such that $\| A\setminus (V_1\cup ...
\cup V_n\cup V) \|_{{\bar \mu },u} <\epsilon $, hence $A\setminus
(V_1\cup ... \cup V_n\cup V)\subset G \setminus G_{\epsilon ,u}$.
Since $V\subset G \setminus G_{\epsilon ,u}$, then $A\cap
G_{\epsilon ,u}\subset V_1\cup ... \cup V_n$. Thus $A\cap
G_{\epsilon ,u}$ is ${\cal R}_{\mu }$-compact.
\par $(2)$. Each function $\phi _u$ is $\cal R$-upper
semi-continuous and $\phi _u$ is bounded from above on each $\cal
R$-compact set, also $ \| A \|_{\mu ,u}\le \| Ch_A \|_{\phi _u}$ for
each $A\in \cal R$. Thus for each $u\in \cal S$ and every $A\in \cal
R$ there is the inequality $\| A\|_{\mu ,u}<\infty $, consequently,
37$(ii)$ is satisfied. \par Take $\epsilon >0$ and a shrinking
subfamily $\cal A$ in $\cal R$ such that $\bigcap_{E\in \cal A} =
\emptyset $. Then the sets $ \{ x\in A: \phi _u(x)\ge \epsilon \} $
form a family $\cal E$ of $\cal R$-compact sets closed under finite
intersections and $\bigcap_{B\in \cal E} B=\emptyset $. Therefore,
there exists $E\in \cal A$ such that $ \{ x\in E: \phi _u(x)\ge
\epsilon \} = \emptyset $. Hence $E\subset \{ x\in G: \phi
_u(x)<\epsilon \} $ and inevitably $ \| A \|_{\mu ,u} <\epsilon $.

\par {\bf 52. Corollary.} {\it If $\mu $ is a $X$-valued measure on
$\cal R$ then for every $\epsilon >0$ and each semi-norm $u\in \cal
S$ the set $\{ x\in G: N_{\mu ,u}(x)\ge \epsilon \} $ is $\cal
R$-locally compact.}

\par {\bf 53. Theorem.} {\it  Let $\mu $ be a measure on $\cal R$
and let also $\cal U$ be a separating covering ring of $G$ being a
sub-ring of ${\cal R}_{\mu }$ and let $\nu $ be a restriction of
$\mu $ on $\cal U$. Then ${\cal U}_{\nu } = {\cal R}_{\mu }$ and
${\bar \nu } = {\bar \mu }$.}
\par {\bf Proof.} Let $u\in \cal S$ be a semi-norm in $X$.
At first we prove that $N_{\mu ,u}(x)\ge N_{\nu ,u}(x)$ for each
$x\in G$. Suppose the contrary, that there exists $y\in G$ such that
$N_{\mu ,u}(y)< N_{\nu ,u}(y)$. This implies that there exists $V\in
\cal R$ such that $\| V \|_{\mu ,u} <N_{\nu ,u}(y)$. \par In
accordance with Lemma 50 $ \| V \|_{\mu ,u} = \| V \|_{{\bar \mu
},u}$. Then for every $x\in G\setminus V$ there exists $B\in \cal U$
such that $y\in B$ and $x\notin B$, since $\cal U$ is the separating
covering ring. Therefore, $ \{ B\setminus V: B\in {\cal U}, y\in B
\} $ is a shrinking sub-collection of ${\cal R}_{\mu }$ whose
intersection is empty. In view of 48$(iii')$ there exists $B\in \cal
U$ such that $y\in B$ and $ \| B\setminus V \|_{{\bar \mu },u}
<N_{\nu ,u}(y)$. Since $ \| V \|_{{\bar \mu },u} <N_{\mu ,u}(y)$.
But $\nu $ is the restriction of $\bar \mu $ on $\cal U$ and $B\in
\cal U$, then $\| B \|_{\nu ,u} \le \| B\|_{{\bar \mu },u} <N_{\nu
,u} (y)$ giving the contradiction, since $y\in B$.
\par Demonstrate now, that ${\cal R}_{\mu }\subset {\cal U}_{\nu }$.
Let $A\in {\cal R}_{\mu }$ and $\epsilon >0$. In view of Lemma 47 it
is sufficient for each $u\in \cal S$ to construct $B\in \cal U$ such
that $N_{\nu ,u}(x)<\epsilon $ for every $x\in A\bigtriangleup B$.
Take $W := \{ x\in A: N_{\mu ,u}(x)\ge \epsilon \} $. In view of
Theorem 51 the set $W$ is ${\cal R}_{\mu }$-compact, consequently,
$\cal U$-compact. Mention that for each $x\in G\setminus W$ there
exists $B\in \cal U$ such that $W\subset B$ with $x\notin B$. Then
the shrinking sub-family $ \{ B\setminus A: B\in {\cal U}, B\supset
W \} $ of ${\cal R}_{\mu }$ has the empty intersection and there
exists $B\in \cal U$ such that $B\supset W$ for which $\| B\setminus
A \|_{{\bar \mu },u}<\epsilon $. Thus $N_{{\bar \mu },u}=N_{\mu ,u}<
\epsilon $ on $B\setminus A$. On the other hand, $N_{\mu
,u}(x)<\epsilon $ on $A\setminus W$ hence on $A\bigtriangleup B$ as
well.
\par Next we show that $\mu $ can be obtained as the restriction of
$\bar \nu $. For $A$, $B$, $\epsilon $, $u$ as in the preceding
paragraph we have $u({\bar \nu }(A) - {\bar \mu }(A)) =u({\bar \nu }
(A) - \nu (B) +{\bar \mu }(B) - {\bar \nu }(A)) = u({\bar \nu
}(A\setminus (A\cap B)) - {\bar \nu } (B\setminus (A\cap B)) + {\bar
\mu } ((B\setminus (A\cap B)) - {\bar \mu }(B\setminus (A\cap B)))
\le \max ( \| A\bigtriangleup B \|_{{\bar \nu},u}, \|
A\bigtriangleup B \|_{{\bar \mu },u})= \| A\bigtriangleup B \|
_{{\bar \mu},u} =\sup_{x\in A\bigtriangleup B} N_{\mu ,u}(x)$.
Therefore, ${\bar \nu }={\bar \mu }$ on ${\cal R}_{\mu }$.
\par Thus ${\cal R}\subset {\cal U}_{\nu }$ and $\mu $ is the
restriction of $\bar \nu $ on $\cal R$. Symmetrically interchanging
$\mu $ and $\nu $ in the proof above one obtains ${\cal U}_{\nu } =
{\cal R}_{\mu }$ and ${\bar \mu } ={\bar \nu }$.

\par {\bf 54. Lemma.} {\it  Let $\mu $ be a $X$-valued measure on
$\cal R$, for a semi-norm $u$ in $X$ and for $\epsilon >0$ put
$G_{\epsilon ,u} := \{ x\in G: N_{\mu ,u}(x)\ge \epsilon \} $. Then
the restriction of the $\cal R$ and ${\cal R}_{\mu }$-topologies to
$G_{\epsilon ,u}$ coincide. A function $f: G\to \sf g$ is ${\cal
R}_{\mu }$-continuous if and only if for each $u\in \cal S$ and
$\epsilon
>0$ the restriction $f|_{G_{\epsilon ,u}}$ is ${\cal R}$-continuous.}
\par {\bf Proof.} In view of Lemma 47 the ${\cal R}$-topology and
${\cal R}_{\mu }$-topology induce the same topology on $G_{\epsilon
,u}$. Therefore, if $f: G\to \sf g$ is ${\cal R}_{\mu }$-continuous,
then it is $\cal R$-continuous on $G_{\epsilon ,u}$.
\par Suppose that $f: G\to \sf g$ has $\cal R$-continuous restrictions
$f|_{G_{\epsilon ,u}}$ for each $u\in \cal S$ and every $\epsilon
>0$. Take any $V$ clopen in $\sf g$. If $A\in {\cal R}_{\mu }$, then
$A\cap G_{\epsilon ,u}$ is $\cal R$-compact due to Theorem 51, as
well as $f^{-1}(V)\cap A\cap G_{\epsilon ,u}$ is $\cal R$-clopen as
a subset of $G_{\epsilon ,u}$. For each $x\in f^{-1}(V)\cap A\cap
G_{\epsilon ,u}$ take $U_x\in \cal R$ with $x\in U_x$ such that
$U_x\cap G_{\epsilon ,u}\subset f^{-1}(V)\cap A\cap G_{\epsilon
,u}$. From this covering of the compact set choose a finite
sub-covering such that $f^{-1}(V)\cap A\cap G_{\epsilon ,u}\subset
U$, where $U := \bigcup_{j=1}^kU_{x_j}$, hence $U\in \cal R$.
Therefore, $f^{-1}(V)\cap A\cap G_{\epsilon ,u}= U\cap G_{\epsilon
,u}$. In view of Lemma 47 $f^{-1}(V)\cap A\in {\cal R}_{\mu }$ for
each $A\in \cal R$. Thus $f^{-1}(V)$ is ${\cal R}_{\mu }$-clopen and
inevitably $f$ is ${\cal R}_{\mu }$-continuous.

\par {\bf 55. Corollary.} {\it If $f: G\to \sf g$ is
${\cal R}_{\mu }$-continuous on each ${\cal R}_{\mu }$-compact set,
then $f$ is ${\cal R}_{\mu }$-continuous on $G$. If $u\in \cal S$
and $E\subset G$ is ${\cal R}_{\mu }$-compact, then $H := \{ x\in E:
N_{\mu ,u}(x)=0 \} $ is finite and there exists $\delta >0$ such
that $N_{\mu ,u}>\delta $ on $E\setminus H$.}
\par {\bf Proof.} In view of Lemma 54 each $\cal R$-compact set
$G_{\epsilon ,u}$ is ${\cal R}_{\mu }$-compact and $f$ is continuous
on every $\cal R$-compact subset of $G_{\epsilon ,u}$. On the other
hand, $G_{\epsilon ,u}$ is $\cal R$-locally compact by Corollary 52.
Therefore, $f$ is $\cal R$-continuous on $G_{\epsilon ,u}$ and
${\cal R}_{\mu }$-continuous on $G$. \par We have that each subset
$A$ of $\{ x\in G: N_{\mu ,u}(x)=0 \} $ is ${\cal R}_{\mu }$-clopen,
since $Ch_A\in L(\mu )$. Take $E$ a ${\cal R}_{\mu }$-compact subset
in $G$. Hence $H$ is finite. Let $\pi \in \bf K$ with $0<|\pi |<1$.
If $\inf \{ N_{\mu ,u}(x): x\in E\setminus H \} =0$, then there
exists a sequence $\{ x_k\in E: k\in {\bf N} \} $ in $E$ such that
$N_{\mu ,u}(x_k)<|\pi |^k$ and $N_{\mu ,u}(x_k) < N_{\mu
,u}(x_{k-1})$ for each $k$. Choose $A_k\in \cal R$ such that $x_k\in
A_k$ and $N_{\mu ,u}(x)<|\pi |^k$ for each $x\in A_k$ and $A_k\cap
A_l=\emptyset $ for each $k\ne l$. Without loss of generality we can
consider the family of non-archimedean semi-norms $\cal S$ in $X$
such that if $u, q\in \cal S$, then $t := \max (u(x),q(x))\in \cal
S$ is also a semi-norm in $X$ (see also \cite{nari}). Therefore, $\|
A \|_{\mu ,t} =\max ( \| A \|_{\mu ,u}, \| A \|_{\mu ,q})$ for each
$A\in \cal R$, as well as $N_{\mu ,t}(x) =\max ( N_{\mu ,u}(x),
N_{\mu ,q}(x))$ for each $x\in G$. Hence $G_{\epsilon ,q}\cap
G_{\epsilon ,u}\supset G_{\epsilon ,t}$ for each $\epsilon >0$ and
$t = \max (u,q)$, $u, q, t\in \cal S$. If $f|_U$ is continuous on a
set $U$ and $W$ is a subset of $U$, $W\subset U$, then evidently
$f|_W$ is continuous.
\par In view of Lemma 54 and Theorem 51 the function \par $g(x) := \sum_k
\pi ^{-k} Ch_{A_k\cap \{ x\in G: N_{\mu ,u}(x)>0 \} } (x)v_k$ \\ is
${\cal R}_{\mu }$-continuous, where $v_k\in X$, $u(v_k)=1$ for each
$k$, since the restriction of $g$ on each $G_{\gamma ,q}$ is
continuous for each $q\in \cal S$ and $\gamma >0$. But $g$ appears
to be not bounded on the ${\cal R}_{\mu }$-compact set $E$. This
gives the contradiction, consequently, $\inf \{ N_{\mu ,u}(x): x\in
E\setminus H \} >0$.

\par {\bf 56. Theorem.} {\it Let $\mu $ be a $X$-valued measure on $\cal R$.
A function $f: G\to \sf g$ is $\mu $-integrable if and only if it
satisfies $(1)$ and $(2)$:
\par $(1)$ $f$ is ${\cal R}_{\mu }$-continuous;
\par $(2)$ for each $u\in \cal S$ and $\epsilon >0$ the set $ \{ x:
x\in G, u(f(x))N_{\mu ,u}(x)\ge \epsilon \} $ is ${\cal R}_{\mu
}$-compact, consequently, contained in some $ \{ x: N_{\mu ,u}(x)\ge
\delta \} $ with $\delta >0$.}
\par {\bf Proof.} If $u$ is a semi-norm in $X$ and $\epsilon >0$,
then the set $G_{\epsilon ,u} := \{ x\in G: N_{\mu ,u}(x)\ge
\epsilon \} $ is $\cal R$-compact by Theorem 51. Without loss of
generality we consider complete $X$ and $\sf g$. For each $f\in
L(G,{\cal R},\mu ,X; {\sf g})$ and each semi-norm $u$ in $X$ there
exists a sequence of simple functions $\{ f_k: k\in {\bf N} \} $
such that $\lim_{k\to \infty } \| f-f_k \|_{\mu ,u} =0$. Then each
$f_k$ is $\cal R$-continuous and the sequence $\{ f_k: k \} $
converges uniformly on $G_{\epsilon ,u}$ to $f$, hence $f$ is $\cal
R$-continuous on $G_{\epsilon ,u}$. In view of Corollary 55 the
function $f$ is ${\cal R}_{\mu }$-continuous.
\par Take a step function $g$ such that $\| f-g \|_{\mu ,u}
<\epsilon $, consequently, $\{ x: u(f(x)) N_{\mu ,u}(x)\ge \epsilon
\} = \{ x: u(g(x)) N_{\mu ,u}(x)\ge \epsilon \} $ and this set is
compact by Theorem 51. Thus from $f \in L(G,{\cal R},\mu ,X; {\sf
g})$ Properties $(1,2)$ follow.
\par Let now Properties $(1,2)$ be satisfied for $f: G\to X$.
For $\delta >0$ and a semi-norm $u$ in $X$ take a ${\cal R}_{\mu
}$-step function $g$ such that $ \| f-g \|_{\mu ,u} < \delta $.
Consider the set $ V := \{ x\in G: u(f(x))N_{\mu ,u}(x)\ge \delta \}
$. The function $N_{\mu ,u}$ is ${\cal R}_{\mu }$-upper
semi-continuous by Theorem 51 and $\sup_{x\in C} N_{\mu ,u}(x) =: w
<\infty $. Since $V$ is compact, then there exists a finite clopen
covering $B_1,...,B_k$ of $V$ such that $u(f(x)-f(y))w<\delta $ for
each $x, y\in B_j\cap V$ with the same $j$, $j=1,...,k$. Choose now
$b_j\in B_j$, hence the set $ \{ x\in G: u(f(x)-f(b_j))N_{\mu
,u}(x)< \delta \} $ is ${\cal R}_{\mu ,u }$-open and contains $B_j$.
Therefore, there exist disjoint sets $W_1,...,W_k$ such that
$W_j\subset \{ x\in G: u(f(x)- f(b_j)) N_{\mu ,u}(x) <\delta \} $
and $B_j := W_j\cap V$, since $G$ is $\cal R$ totally disconnected
with the clopen base of its topology. Take the step function $g (x)
:= \sum_{j=1}^k f(b_j) Ch_{W_j}(x)$. For $x\in W_j$ we have
$u(f(x)-g(x)) N_{\mu ,u}(x) = u(f(x) - f(b_j)) N_{\mu ,u}(x)< \delta
$. At the same time for $x\notin  \bigcup_{j=1}^k W_j$ we have
$u(f(x) -g(x)) N_{\mu ,u}(x) = u(f)) N_{\mu ,u}(x) <\delta $,
consequently, $ \| f-g \|_{\mu ,u} \le \delta $.

\par {\bf 57. Corollary.} {\it Let $\mu $ be a $X$-valued measure on
$\cal R$, $g\in L(G, {\cal R}, \mu ,X; {\sf g})$, $f: G\to {\sf g}$
be ${\cal R}_{\mu }$-continuous and $u(f(x)) \le u(g(x))$ for each
semi-norm $u\in {\sf g}$ and for every $x\in G$, then $f\in L(G,
{\cal R}, \mu ,X; {\sf g})$.}

\par {\bf 58. Corollary.} {\it The space
$L(G,{\cal R},\mu ,X; {\sf g})$ is complete and locally $\bf
K$-convex. If $X$ and $\sf g$ are normed spaces, then $L(G,{\cal
R},\mu ,X; {\sf g})$ is the Banach space.}
\par {\bf Proof.} By the construction above $L(G,{\cal R},\mu ,X; {\sf g})$
is the completion of the space of step functions relative to the
family of semi-norms $ \| * \|_{\mu ,u}$, where $u$ is a semi-norm
in $X$, $\sf g$. Therefore, $L(G,{\cal R},\mu ,X; {\sf g})$ is
isomorphic with $L(G,{\cal R},\mu ,{\tilde X}; {\sf g})$ and
complete, where $\tilde X$ is the completion of $X$ and ${\tilde
{\sf g}}$ is the completion of $\sf g$ as the $\bf K$-convex spaces.
Particularly, when $X$ and $\sf g$ are normed spaces, then $\tilde
X$ and ${\tilde {\sf g}}$ and $L(G, {\cal R}, \mu ,{\tilde X};
{\tilde {\sf g}})$ are the Banach spaces.

\par {\bf 59. Definition.} If $X$ is a normed space over $\bf K$
and $x_1, x_2,...$ is a sequence in $X$ such that $\| a_1x_1 + ... +
a_nx_n \| \ge t \max \{ \| a_jx_j \| : j=1,...,n \} $ for each
$a_1,...,a_n\in \bf K$ and $n\in \bf N$ not exceeding the length of
the sequence, where $0<t\le 1$ is a marked number, then $\{
x_1,x_2,... \} $ are called $t$-orthogonal. If $t=1$, then $\{ x_1,
x_2,... \} $ are called orthogonal.
\par Naturally in the case $t=1$ the inequality, $\ge $, reduces to the
equality, $=$, due to the non-archimedean property of the norm.

\par {\bf 60. Theorem.} {\it Let $\mu $ be a $X$-valued measure on
$(G,{\cal R})$ and let a Banach algebra $\sf g$ has a
$t_0$-orthogonal basis, where $0<t_0\le 1$. Then $L(G,{\cal R},\mu
,X;{\sf g})$ has a $t$-orthogonal basis with $0<t<t_0$. If the
valuation group of the field $\bf K$ is discrete in $(0,\infty )$,
then $L(G,{\cal R},\mu ,X;{\sf g})$ has an orthogonal basis.}
\par {\bf Proof.} If the valuation group of $\bf K$ is discrete,
then  the Banach space over $\bf K$ has an orthogonal basis by
Theorem 5.16 \cite{roo}, particularly, for $L(G,{\cal R},\mu ,X;{\sf
g})$ due to Corollary 58.
\par In general, suppose that the valuation group of $\bf K$ is
dense in $(0,\infty )$ and $N_{\mu }(x)>0$ for each $x\in G$. Take a
marked $0<t<t_0$ such that $\sf g$ has the $t_0$-orthogonal basis.
Then choose $\pi \in \bf K$ such that $t_1<|\pi |<1$, where
$t_1=t/t_0$, and define the function $h: G\to \bf K$ such that
$h(x)=\pi ^n$, when $|\pi |^{n+1}<N_{\mu }(x)\le |\pi |^n$ and $n\in
\bf Z$, consequently, $|\pi | |h(x)| < N_{\mu }(x)\le |h(x)|$ for
each $x\in G$. Denote by $G_d$ the set $G$ in the discrete topology
and define the mapping $q: L(G,{\cal R},\mu ,X;{\sf g})\ni f\mapsto
hf\in BC(G_d,{\sf g})$, where $BC(Y,W)$ denotes the space of all
bounded continuous mappings from a topological space $Y$ into a $\bf
K$-linear normed space $W$. The space $BC(G_d,{\sf g})$ is supplied
with the norm $\| v \|_{\infty } := \sup_{x\in G_d} \| v(x) \| $,
where $\| * \|$ is the norm in $\sf g$. Therefore, \par $(1)$ $|\pi
| \| qf \|_{\infty }< \| f \|_{\mu } \le \| qf \|_{\infty }$ \\ for
each $f\in L(G,{\cal R},\mu ,X;{\sf g})$. If $A\in \cal R$ and $b\in
\sf g$, then $hbCh_A\in BC(G_d,{\sf g})$, since $\| h(x)bCh_A(x):
x\in G \} \subset \{ \pi ^n: n\in {\bf Z}; |\pi |^{n+1}< \| b \|
|Ch_A(x)| \} $. \par By $\bf K$ linearity and Property $(1)$ of the
mapping $q$ the range $q(L(G,{\cal R},\mu ,X;{\sf g}))$ is the
closed $\bf K$-linear subspace in $BC(G_d,{\sf g})$. In the space
$BC(G_d,{\sf g})$ the product $BC(G_d,{\bf K})\times {\sf g}$ is
everywhere dense, while $L(G,{\cal R},\mu ,X;{\bf K})\times {\sf g}$
is everywhere dense in $L(G,{\cal R},\mu ,X;{\sf g})$. In view of
Corollaries 5.23 and 5.25 \cite{roo} $BC(G_d,{\bf K})$ has an
orthogonal basis, hence $BC(G_d,{\sf g})$ has a $t_0$-orthogonal
basis. \par In accordance with the Gruson Theorem 5.9 \cite{roo} if
$E$ is a Banach space with an orthogonal basis, then each its closed
$\bf K$-linear subspace has an orthogonal basis. The space
$q(L(G,{\cal R},\mu ,X;{\bf K}))$ is closed in $BC(G_d,{\bf K})$,
hence it has an orthogonal basis $\{ e_j: j \} $. Thus $L(G,{\cal
R},\mu ,X;{\sf g})$ has the $t$-orthogonal basis $q^{-1}(e_j\times
s_k)$, where $\{ s_k: k\} $ is the $t_0$-orthogonal basis in $\sf
g$.

\par {\bf 61. Theorems.} {\it Let $X$ and $\sf g$ be as in \S 39
and in addition let $X$ be an algebra over $\bf K$ with a family of
multiplicative semi-norms. Suppose that $\mu $ and $\nu $ are
$X$-valued measures on separating covering rings $\cal R$ of a set
$G$ and $\cal T$ of a set $H$. Then
\par $(1)$ the finite unions of the sets $A\times B$, $A\in \cal R$,
$B\in \cal T$ form the separating covering ring ${\cal R}\otimes
{\cal T}$ of $G\times H$; \par $(2)$ there exists a unique measure
$\mu \times \nu $ on ${\cal R}\times {\cal T}$ such that $\mu \times
\nu (A\times B)=\mu (A)\nu (B)$ for each $A\in \cal R$ and $B\in
\cal T$, $N_{\mu \times \nu ,u}(x,y)=N_{\mu ,u}(x) N_{\nu ,u}(y)$
for each $x\in G$, $y\in H$ and every semi-norm $u\in \cal S$ in
$X$; \par $(3)$ if $f\in L(G\times H, {\cal R}\times {\cal T},\mu
\times \nu , X; {\sf g})$, then $H\ni y\mapsto \int_G f(x,y)\mu
(dx)$ is a $\nu $-almost everywhere defined $\nu $-integrable
function on $H$ and $G\ni x\mapsto \int_H f(x,y)\nu (dy)$ is a $\mu
$-almost everywhere defined $\mu $-integrable function on $G$ and
\par $\int_{G\times H}f(x,y)\mu \times \nu (dx,dy)= \int_H(\int_Gf(x,y)
\mu (dx))\nu (dy)$. Moreover, if $X$ is commutative, then
$\int_H(\int_Gf(x,y) \mu (dx))\nu (dy) = \int_G(\int_Hf(x,y) \nu
(dx))\mu (dy)$;
\par $(4)$ in particular, if $\sf g$ and $X$ are commutative,
$f\in L(G,{\cal R},\mu ,X; {\sf g})$ and $h\in L(H,{\cal T},\nu ,X;
{\sf g})$, then $f(x)h(y)\in L(G\times H, {\cal R}\times {\cal
T},\mu \times \nu , X; {\sf g})$ and \par $\int_{G\times H} f(x)g(y)
\mu (dx)\nu (dy) = (\int_G f(x)\mu (dx)) (\int_H g(y)\nu (dy))$;
\par $(5)$ $L(G\times H, {\cal R}\times {\cal T},\mu
\times \nu , X; {\sf g})$ is $\bf K$-linearly topologically
isomorphic with the tensor product $L(G, {\cal R},\mu , X; {\sf g})
{\hat \otimes } L(H,{\cal T},\nu , X; {\sf g})$.}

\par {\bf Proof.} $(1)$. If $\cal U$ and $\cal V$ are coverings
of $G$ and $H$ by elements from $\cal R$ and $\cal T$ respectively,
then $\bigcup_{A\in {\cal U}, B\in {\cal V}} A\times B= G\times H$.
If $(x_1,y_1)\ne (x_2,y_2)\in G\times H$, then either $x_1\ne x_2$
or $y_1\ne y_2$. In the first case take $A\in \cal R$ such that
$x_1\in A$ and $x_2\in G\setminus A$ and $x_1\in B$ with $B\in \cal
T$, then $A\times B$ separates them: $(x_1,y_1)\in A\times B$ and
$(x_2,y_2)\notin A\times B$.
\par $(2)$. Put $\mu \times \nu (E) := \int_{G\times H} Ch_C(x,y)
\mu (dx)\nu (dy)$ for each $E\in {\cal R}\times {\cal T}$, hence
$\mu \times \nu $ is additive. For $N_u(x,y) := N_{\mu ,u}(x) N_{\nu
,u}(y)$ and each $A\in \cal R$, $B\in \cal T$ we get $u((\mu \times
\nu )) (A\times B))\le \| Ch_{A\times B} \|_{N_u}$, consequently, $
\| C \|_{\mu \times \nu ,u} \le \| Ch_C \|_{N_u}$. \par Naturally
$G$ and $H$ and $G\times H$ are supplied with the $\cal R$ and $\cal
T$ and ${\cal R}\times \cal T$ topologies. Then $ \{ (x,y)\in
A\times B: N_u(x,y)\ge \epsilon \} \subset \{ x\in A: N_{\mu
,u}(x)\ge \epsilon \} \times \{ y\in B: N_{\nu ,u}(y)\ge \epsilon \}
$ for each $A\in \cal R$ and $B\in \cal T$ and $\epsilon >0$. The
function $N_u$ is upper semi-continuous and $\{ x\in A: N_{\mu
,u}(x)\ge \epsilon \} \times \{ y\in B: N_{\nu ,u}(y)\ge \epsilon \}
$ is compact, consequently, $ \{ (x,y)\in A\times B: N_u(x,y)\ge
\epsilon \}$ is compact and inevitably $\{ (x,y)\in C: N_u(x,y)\ge
\epsilon \} $ is compact for each $C\in {\cal R}\times \cal T$ and
every $\epsilon >0$. Therefore, by Theorem 51 $\mu \times \nu $ is
the measure and $N_{\mu \times \nu ,u}(x,y)\le N_u(x,y)$ for each
$u\in \cal S$ and each $x\in G$ and every $y\in H$. Each $u$ is the
multiplicative semi-norm in $X$ and $u((\mu \times \nu )(A\times B))
= u(\mu (A)) u(\nu (B))$ and hence $\sup \{ u((\mu \times \nu )(C)):
C\in {\cal R}\times {\cal T}, C\subset A\times B \} =  \| A\times B
\| _{\mu \times \nu ,u} \ge \| A \|_{\mu ,u} \| B \|_{\nu ,u}$ for
each $A\in \cal R$ and every $B\in \cal T$. Therefore, $N_{\mu
\times \nu ,u}(x,y)\ge N_u(x,y)$ for each $u$, $x$ and $y$.
\par $(3)$. If $f$ is a step function, then $(3)$ is satisfied
due to $(2)$. If $f\in L(G,{\cal R},\mu ,X; {\sf g})$, then for each
$u$ there exists a sequence of step functions $f_1, f_2,...$
converging to $f$ such that $\| f-f_n \|_{N_u} \le 1/n$ for each
$n\in \bf N$. Since $u(f(x,y)-f_n(x,y))N_{\mu ,u}(x)N_{\nu ,u}(y)\le
1/n$ for each $x, y$. Thus $f(*,y)\in L(G,{\cal R},\mu ,X;{\sf g})$
for $\nu $-almost every $y\in H$, consequently, $u(\int_Gf(x,y)\mu
(dx) -\int_Gf_n(x,y)\mu (dx)) N_{\nu ,u}(y)\le 1/n$. Therefore, the
function $H\ni y\mapsto \int_Gf(x,y)\mu (dx)$ is defined for $\nu
$-almost all $y\in H$ and is $\nu $-integrable. Since \par
$\int_H(\int_G f_n(x,y)\mu (dx))\nu (dy) =\int_{G\times
H}f_n(x,y)\mu \times \nu (dx,dy)$ for each $n$, then
\par $u(\int_H(\int_G f_n(x,y)\mu (dx))\nu (dy)- \int_{G\times
H}f_n(x,y)\mu \times \nu (dx,dy))\le 1/n $ for each $n\in \bf N$ and
each $u\in \cal S$, consequently, \par $\int_H(\int_G f(x,y)\mu
(dx))\nu (dy) =\int_{G\times H}f(x,y)\mu \times \nu (dx,dy)$.
\par Part $(4)$ follows from $(2,3)$ as the particular case.
\par $(5)$. There is a bilinear continuous mapping ${\sf m}:
L(G, {\cal R},\mu , X; {\sf g}) \times L(H,{\cal T},\nu , X; {\sf
g})\ni (f,h)\mapsto fh\in L(G\times H, {\cal R}\times {\cal T},\mu
\times \nu , X; {\sf g})$. If $X$ is a Banach space, then the norm
of $\sf m$ is $\| {\sf m} \| =1$. \par Let $Y$ be a complete locally
$\bf K$-convex space over $\bf K$ and let $F: L(G, {\cal R},\mu , X;
{\sf g}) \times L(H,{\cal T},\nu , X; {\sf g})\to Y$ be a continuous
$\bf K$-bilinear mapping. Then there exists the mapping $F_{\sf
m}(Ch_{A\times B})=F_{\sf m}(Ch_A \otimes Ch_B)= F(Ch_A,Ch_B)\in Y$
for each $A\in \cal R$ and $B\in \cal T$. A step function $f:
G\times H\to {\sf g}$ write in the form $f(x,y)=\sum_j b_j
Ch_{A_j\times B_j}(x,y)$, where $b_j\in \sf g$ and $A_j\in \cal R$,
$B_j\in \cal T$ for each $j$, $(A_j\times B_j)\cap (A_i\times B_i) =
\emptyset $ for each $i\ne j$. \par For each semi-norm $v$ in $Y$
and each semi-norm $u$ in $X$ there exists $k$ such that
\par $\sup_{x\in G, y\in H} v(F_{\sf m}(b_kCh_{A_k\times B_k})(x,y))
N_{\mu \times \nu ,u}(x,y)$ \\  $ = \max_j v(F_{\sf
m}(b_jCh_{A_j\times B_j})(x,y)) N_{\mu \times \nu ,u}(x,y)$, hence
\par $v(F_{\sf m}f (x,y))N_{\mu \times \nu ,u}(x,y)\le \sup_{x\in G,
y\in H} v(F_{\sf m}(b_kCh_{A_k\times B_k})(x,y)) N_{\mu \times \nu
,u}(x,y)\le u(b_j) \| F \|_{u,v} \| Ch_{A_k} \|_{\mu ,u} \| Ch_{B_k}
\| _{\nu ,u}\le \| F \|_{u,v} \| f \|_{\mu \times \nu , u}$,\\ where
$ \| F \|_{u,v} := \sup_{w\in X, u(w)>0} v(Fw)/u(w)$. \par
Therefore, the mapping $F_{\sf m}$ has the continuous extension
$F_{\sf m}$ from $L(G\times H, {\cal R}\times {\cal T},\mu \times
\nu , X; {\sf g})$ into $Y$. Thus $F_{\sf m}(f\otimes h)=F(f,h)$ for
each $f\in L(G, {\cal R},\mu ,X; {\sf g})$ and $h\in L(H,{\cal
T},\nu , X; {\sf g})$ and inevitably $L(G\times H, {\cal R}\times
{\cal T},\mu \times \nu , X; {\sf g})$ is $\bf K$-linearly
topologically isomorphic with the tensor product $L(G, {\cal R},\mu
, X; {\sf g}) {\hat \otimes } L(H,{\cal T},\nu , X; {\sf g})$.

\par {\bf 62. Theorem.} {\it  Let ${\sf g}_j$ and $X_j$ be a family of
algebras and locally convex spaces over $\bf K$ satisfying
Conditions of \S 39, $j\in \beta $, where $\beta $ is a set. Suppose
that for each $j$ there is a measure $\mu _j: {\cal R}\to X_j$,
$X=\bigotimes_{j\in \beta } X_j$ and ${\sf g} = \bigotimes_{j\in
\beta } {\sf g}_j$ are supplied with the product topologies. Then
there exists a measure $\mu = \bigotimes_{j\in \beta } \mu _j$ on
${\cal R}$ with values in $X$ such that $L(G,{\cal R},\mu ,X;{\sf
g})$ is the completion of the direct sum $\bigoplus_{j\in \beta }
L(G,{\cal R},\mu _j, X_j; {\sf g}_j)$.}
\par {\bf Proof.} Since each $X_j$ and every ${\sf g}_j$ are
complete, then $X$ and $\sf g$ are complete (see \cite{eng,nari}).
Naturally $X$ is the locally $\bf K$-convex space and $\sf g$ is the
algebra such that $x+y = (x_j+y_j: j)$ and $ab = (a_jb_j: j)$, where
$x, y\in X$, $a, b\in \sf g$, $x = (x_j: j)$, $a = (a_j: j)$,
$ax=(a_jx_j: j)$ (see Theorem 5.6.1 \cite{nari}). Thus $X$ is the
unital left $\sf g$-module. For each $j\in \beta $ there are defined
the projectors $\pi _j(x)=x_j$ and $\pi _j(a)=a_j$ on $X$ and $\sf
g$. Therefore, topologies of $X$ and $\sf g$ are characterized by
the families of semi-norms $u$ such that $u(x)=\max \{ u_j(x_j):
j\in \alpha \} $ and $u(b)=\max \{ u_j(b_j): j\in \alpha \} $, where
$\alpha $ is a finite subset in $\beta $, $u_j\in {\cal S}_j$, where
${\cal S}_j$ is a consistent family of semi-norms $u_j$ in $X_j$ and
${\sf g}_j$ denoted by the same symbol for shortening the notation.
\par If $A, B\in \cal R$ and $A\cap B=\emptyset $, then $\mu (A\cup
B)= (\mu _j(A\cup B): j) =(\mu _j(A) +\mu _j(B): j)= (\mu _j(A): j)
+ (\mu _j(B): j)= \mu (A) +\mu (B)$, hence $\mu $ is additive. If
$u$ is a semi-norm in $X$ and $A\in \cal R$, then for each $C\subset
A$, $C\in \cal R$ there is the inequality $u(\mu (C)) = \max \{
u_j(\mu _j(C)): j\in \alpha \} <\infty $, since $\alpha $ is a
finite set and each $\mu _j$ is bounded. If $\cal A$ is a shrinking
family in $\cal R$ and $\bigcap_{A\in \cal A} A =\emptyset $, then
$\lim_{A\in {\cal A}}u_j(\mu _j(A))=0$ for each $j$ and each $u_j\in
{\cal S}_j$, hence $\lim_{A\in {\cal A}}u(\mu (A))=0$ for each
semi-norm $u\in \cal S$ in $X$. Thus $\mu $ is the measure on $\cal
R$ with values in $X$.
\par If $f\in L(G,{\cal R},\mu ,X; {\sf g})$, then for each $u\in
\cal S$ there exists a sequence $\{ f_n: n\in {\bf N} \} $ of simple
functions such that $\| f-f_n \|_{\mu ,u} \le 1/n$. We have that
$\bigoplus_j{\sf g}_j$ is everywhere dense in $\sf g$, where
elements of the direct sum as usually are $b= (b_j: j\in {\beta },
b_j\in {\sf g}_j \} $ such that the set $ \{ j: b_j\ne 0 \} $ is
finite (see Example 5.10.6 in \cite{nari}). Thus each simple
function can be chosen taking values in $\bigoplus_{j\in \beta }
{\sf g_j}$. Therefore, $L(G,{\cal R},\mu ,X;{\sf g})$ is the
completion of the direct sum $\bigoplus_{j\in \beta } L(G,{\cal
R},\mu _j, X_j; {\sf g}_j)$.
\par Certainly if $\beta $ is finite, then the direct sum and
the direct product coincide.
\par {\bf 63. Corollary.} {\it If suppositions of Theorem 62 are
satisfied and $f\in L(G,{\cal R},\mu ,X; {\sf g})$, then
$\int_Gf(x)\mu (dx) = (\int_G f_j(x)\mu _j(dx_j): j\in \beta )$,
where $f_j\in L(G,{\cal R},\mu ,X_j; {\sf g}_j)$ for each $j$.}
\par {\bf Proof.} The formula of this Corollary is satisfied for each $f$
in \\ $\bigoplus_{j\in \beta } L(G,{\cal R},\mu _j, X_j; {\sf
g}_j)$, but the latter space is everywhere dense in $L(G,{\cal
R},\mu ,X; {\sf g})$ by Theorem 61. The mapping $\int_G : L(G,{\cal
R},\mu ,X; {\sf g})\to X$ is continuous: $u(\int_Gf(x)\mu (dx))\le
\| f \|_{\mu ,u}$ due to Lemma 45, consequently, the statement of
this Corollary follows by the continuity.

\par {\bf 64. Theorem.} {\it Let $\sf g$, $X$, $G$, $\cal R$, $\mu $
be as in \S 39 and let $F$ be a continuous homomorphism of a left
unital $\sf g$-module $X$ into a uniformly complete left unital $\sf
h$-module $Y$. Then $F$ induces a continuous homomorphism ${\hat F}:
L(G,{\cal R},\mu ,X; {\sf g})\to L(G,{\cal R},\nu ,Y; {\sf h})$ such
that $F(\int_Gf(x)\mu (dx)) = \int_G {\hat F}(f)\nu (dx)$ for each
$f\in L(G,{\cal R},\mu ,X; {\sf g})$, where $\nu =F(\mu )$ is an
$Y$-valued measure. If $F({\sf g})=\sf h$, then $\hat F$ is
epimorphic.}

\par {\bf Proof.} If $A, B\in \cal R$ with $A\cap B=\emptyset $,
then $\nu (A\cup B) := F(\mu (A\cup B)) = F(\mu (A)) + F(\mu
(B))=\nu (A) +\nu (B)\in Y$. If $A\in \cal R$ and $u$ is a semi-norm
in $X$, then $\sup_{C\subset A, C\in \cal R} u(\mu (C))<\infty $,
hence for each semi-norm $v$ in $Y$ there is the inequality
$\sup_{C\subset A, C\in \cal R} v(\nu (C))<\infty $, since $F$ is
continuous. If $\cal A$ is a shrinking family in $\cal R$ with
$\bigcap_{A\in \cal A} A= \emptyset $, then $\lim_{\cal A} \mu
(A)=0$ hence $0=\lim_{A\in \cal A} F(\mu (A))=\lim_{A\in \cal A} \nu
(A)$, since $F$ is continuous. Thus $\nu $ is the $Y$-valued
measure.
\par If $v$ is a semi-norm in $Y$, then $v_F(q):= v(F(q))$ is the
continuous semi-norm in $X$, where $q\in X$, hence $v_F(\mu (A)) =
v(\nu (A))$ for each $A\in \cal R$ and inevitably $N_{\mu ,v_F}(x)=
N_{\nu ,v}(x)$ for each $x\in G$.
\par If $f: G\to \sf g$ is a step function $f(x) =\sum_j a_j
Ch_{A_j}(x)$, where $a_j\in \sf g$, $A_j\in \cal R$, then ${\hat
F}(f) = F(f)= \sum_j F(a_j) Ch_{A_j}(x)$, since $F(a_1b_1+a_2b_2)=
F(a_1)F(b_1)+ F(a_2)F(b_2)$ for each $a_j, b_j\in \sf g$, also $0,
1\in \sf g$, $F(0)=0$, $F(1)=1\in \sf h$. Moreover, $\| F(f) \|_{\nu
,v} = \| f \|_{\mu ,v_F}$ for each step function $f$ and each
semi-norm $v$ in $Y$, consequently, $\hat F$ is continuous and $\bf
K$-linear and has the continuous extension ${\hat F}: L(G,{\cal
R},\mu ,X; {\sf g})\to L(G,{\cal R},\nu ,Y; {\sf h})$. \par For each
$s\in L(G,{\cal R},\nu ,Y; {\sf h})$ and each semi-norm $v$ in $Y$
take a sequence $s_n$ of simple functions $s_n: G\to \sf h$
converging to $s$ such that $ \| s-s_n \| _{\nu ,v}\le 1/n$ for each
$n$. If $F({\sf g})=\sf h$, then for each $s_n=\sum_j b_{j,n}
Ch_{A_j}$ there exists a simple function $f_n=\sum_j a_{j,n}
Ch_{A_j}$ such that $f_n: G\to \sf g$ and ${\hat F}(f_n)=s_n$,
$a_{j,n}\in \sf g$ for each $j, n$, where $b_{j,n}\in \sf h$,
$A_j\in \cal R$. But $\{ f_n: n \} $ is a fundamental sequence
relative to $\| * \|_{\mu ,v_f}$, consequently, there exists $f\in
L(G,{\cal R},\mu ,X; {\sf g})$ such that ${\hat F}(f)=s$.
\par  By the conditions of this theorem $F(a_1w_1+a_2w_2) = F(a_1)F(w_1)
+ F(a_2) F(w_2)$ for each $a_1, a_2\in \sf g$, $w_1, w_2\in X$,
where $F(a_j)\in \sf h$ and $F(w_j)\in Y$. Therefore, for each step
function $f(x) = \sum_ja_jCh_{A_j}(x)$ we have $F(\int_G f(x)\mu
(dx)) = F(\sum_j a_j\mu (A_j)) = \sum_j F(a_j)\nu (A_j) = \int_G
{\hat F}(f)(x)\nu (dx)$ and $v_F(\int_Gf(x)\mu (dx)) = v(\int_G
{\hat F}(f)(x)\nu (dx))$ for each semi-norm $v$ in $Y$,
consequently, by the continuity $F(\int_G f(x)\mu (dx)) = \int_G
{\hat F}(f)(x)\nu (dx)$ for each $f\in L(G,{\cal R},\mu ,X; {\sf
g})$.

\par {\bf 65. Definition.} Let ${\cal R}={\sf Bco}(G)$ be the ring
of all clopen subsets of a zero-dimensional Hausdorff space. A
measure $\mu : {\sf Bco}(G)\to X$ is called a tight measure, where
$X$ is as in \S 39. The family $M=M(G,X)$ of all such tight measures
form the $\bf K$-linear space with the family of semi-norms
\par $\| \mu \|_u = \sup_{A\in {\sf Bco}(G)} u(\mu (A)) = \| G \|_{\mu
,u} = \sup_{x\in G} N_{\mu ,u}(x)$, \\
where $u\in \cal S$ is a semi-norm in $X$. In particular, if $X$ is
the normed space, then $M(G,X)$ is the normed space.
\par  The closure of the set $\{ x\in G:
\exists u\in {\cal S}, N_{\mu ,u}(x)>0 \} $ we call the support of
the measure $\mu $.

\par {\bf 66. Theorem.} {\it If $\cal R$ is a covering ring of $G$
being the base of the zero-dimensional Hausdorff topology in $G$ and
if $\mu $ is a $X$-valued measure on $\cal R$, then $(f\mu )(A) :=
\int_G Ch_A(x)f(x)\mu (dx)$ is the tight measure for $f\in L(G,{\cal
R},\mu ,X;{\bf K})$ and the mapping $\psi _{\mu } := \psi :
L(G,{\cal R},\mu ,X;{\bf K})\ni \mapsto (f\mu )$ is the $\bf
K$-linear topological embedding $L(G,{\cal R},\mu ,X;{\bf K})$ into
$M(G,X)$.}
\par {\bf Proof.} For each $u\in \cal S$ the set $\{ x\in G:
N_{\mu ,u}(x)>0 \} $ is $\sigma $-compact, that is the countable
union of compact subsets by Theorem 51. Therefore, if $A$ is a
clopen subset in $G$, then $Ch_A\in L(G,{\cal R},\mu ,X;{\bf K})$,
since for each $u\in \cal S$ and $\epsilon >0$ there exists a
sequence $f_n\in L(G,{\cal R},\mu ,X;{\bf K})$ with $\| Ch_A- f_n
\|_{\mu ,u} \le 1/n$ and supports $supp (f_n)\supset \{ x\in G:
N_{\mu ,u}(x)\ge 1/n \} $. Thus $f\mu $ is defined on ${\cal
R}_{f\mu }\supset {\sf Bco}(G)$, consequently, $f\mu \in M(G,X)$.
Evidently $\psi (af+bg) = a\psi (f) + b\psi (g) = af\mu + bg\mu $.
On the other hand, $N_{f\mu ,u}(x)\le \| f \|_{\mu ,u} N_{\mu
,u}(x)$, hence $\psi $ is continuous. In view of Theorem 56 $u(f(x))
N_{\mu ,u}(x)= N_{f\mu ,u}(x)$ for each $u\in \cal S$ and $x\in G$,
hence $\psi $ is the topological embedding. If $X$ is a Banach
space, then $\psi $ is the isometric embedding.

\par {\bf 67.} Let $Y^*$ denote the topological dual space
of all continuous $\bf K$-linear functionals on a $\bf K$-linear
space $Y$, $Mat_n({\bf K})$ denotes the algebra of all square
$n\times n$ matrices, $n\in \bf N$, $BC(G,Y)$ denotes the space of
all continuous bounded functions from $G$ into $Y$.
\par {\bf Theorem.} {\it If $G$ is a zero-dimensional Hausdorff
space and $\mu \in M(G,Mat_n({\bf K}))$, then $BC(G,Mat_n({\bf K}))
\subset L(G,{\cal R},\mu ,Mat_n({\bf K});Mat_n({\bf K}))$ and
$\lambda _{\mu }(f) := \int_G f(x)\mu (dx)$ provides the $\bf
K$-linear isometric embedding $\lambda : M(G,Mat_n({\bf
K}))\hookrightarrow BC (G,Mat_n({\bf K}))^*$. If $G$ is compact,
then $\lambda $ is the isomorphism.}
\par {\bf Proof.} In view of Theorem 56 we get the inclusion
$BC(G,Mat_n({\bf K})) \subset L(G,{\cal R},\mu ,Mat_n({\bf
K});Mat_n({\bf K}))$. Moreover, the mapping $\lambda $ is defined
and it is $\bf K$-linear. \par Since $Mat_n({\bf K})$ is finite
dimensional over $\bf K$, then its topological dual space is
isomorphic with $Mat_n({\bf K})$. It has the natural norm topology
$\| b \| := \max_{1\le i, j \le n} |b_{i,j}|$. If $\mu \in
M(G,Mat_n({\bf K}))$, then $ \| \mu \| = \sup_{A\in {\sf Bco}(G)}
\max_{i,j} |\mu _{i,j}(A)| = \sup_{A\in {\sf Bco}(G)} |\lambda _{\mu
} (Ch_A)| \le \| \lambda _{\mu } \| $, consequently, $\lambda $ is
the isometric embedding.

\par If $q\in BC(G,Mat_n({\bf K}))^*$, then $\mu (A)$ having matrix
elements $q(E_{i,j} Ch_A) =: \mu ^q_{i,j}(A)$ is and additive
function on ${\sf Bco}(G)$ with values in $Mat_n({\bf K})$, where
$E_{i,j}$ is the $n\times n$ matrix with $1$ at the $(i,j)$-th place
and zeros at others places. For each $A\in {\sf Bco}(G)$ and $b\in
Mat_n({\bf K})$ we have $bCh_A\in BC(G,Mat_n({\bf K}))$. Therefore,
$\mu ^q$ is defined on ${\sf Bco}(G)$. If ${\cal A}\subset {\sf
Bco}(G)$ is a shrinking family and $G$ is compact, then from
$\bigcap_{A\in \cal A}A =\emptyset $ it follows, that $\emptyset \in
\cal A$, since each $A\in \cal A$ is closed in $G$. Since $q$ is
continuous, then for compact $G$ the mapping $\mu ^q$ is the
measure. In this particular case $BC(G,Mat_n({\bf K}))$ is
isomorphic with the space $C(G,Mat_n({\bf K}))$ of all continuous
functions from $G$ into $Mat_n({\bf K})$.

\par {\bf 68. Theorem.} {\it A function $f: G\to \sf g$ is $\mu
$-integrable for each $\mu \in M(G,X)$ if and only if $f$ is bounded
and for each compact subset $V$ in $G$ the restriction $f|_V$ of $f$
to $V$ is continuous.}
\par {\bf Proof.} Suppose that $f$ is $\mu
$-integrable for each $\mu \in M(G,X)$. Take a semi-norm $u\in \cal
S$ in $X$ and a number $\pi \in \bf K$ such that $0<|\pi |<1$. If
$f$ is not bounded, then there exists a sequence $b_j\in \bf G$ such
that $b_i\ne b_j$ for each $i\ne j$ and $\lim_{n\to \infty } u(\pi
^n f(b_n)) =\infty $. Put $\mu := \sum_n \pi ^n x_n \delta _{b_n}$,
where $\delta _b(A) := 1$ if $b\in A$, $\delta _b(A) =0$ if $b\notin
A$, $x_n\in X$, $u(x_n)=1$. Therefore, $\mu \in M(G,X)$ and $N_{\mu
,u}(b_n)= |\pi |^n$ for every $n\in \bf N$, hence $\| f \|_{\mu
,u}=\infty $ and $f\notin L(G,{\sf Bco}(G), \mu ,X; {\sf g})$. In
view of Theorem 7.9 \cite{roo} if $V$ is a compact subset in $G$,
then there exists a measure $\lambda : {\sf Bco}(G)\to \bf K$ such
that $N_{\lambda }(x)=1$ for each $x\in V$ and $N_{\lambda }(x)=0$
for each $x\in G\setminus V$. Take $y\in X$ with $u(y)=1$, then $\mu
=y\lambda $ is the $X$-valued measure and $N_{\mu ,u} = Ch_V$.
Therefore, due to Lemma 54 and Theorem 56 $f|_V$ is continuous.
\par Suppose now that $f$ is bounded and its restriction to each
compact subset of $G$ is continuous. Take any $\mu \in M(G,X)$, then
due to Corollary 55 the mapping $f$ is ${\sf Bco}(G)_{\mu
}$-continuous. If $z\in \sf g$, then $q_z\in L(G,{\sf Bco}(G),\mu
,X;{\sf g})$, where $q_z(x)=z$ for every $x\in G$. Take $z\in \sf g$
such that $u(f(x))\le u(z)$ for each $x\in G$, consequently, in view
of Corollary 57 $f\in L(G,{\sf Bco}(G),\mu ,X; {\sf g})$.

\par {\bf 69. Definition.} Let $G$ be a zero-dimensional Hausdorff
topological space. If for each $U\subset G$ it is clopen in $G$ if
and only if $U\cap V$ is clopen in $V$ for each compact subset $V$
in $G$, then $G$ is called the $k_0$-space.

\par {\bf 70. Corollary.} {\it If $G$ is a $k_0$-space, then
$BC(G,{\sf g})= \bigcap_{\mu \in M(G,X)} L(G,{\sf Bco}(G),\mu
,X;{\sf g})$.}
\par {\bf Proof.} By Theorem 3.3.21 \cite{eng} a mapping $f: G\to Y$
of a $k$-space $G$ into a topological space $Y$ is continuous if and
only if for each compact subset $V$ in $G$ the restriction $f|_V$ of
$f$ to $V$ is continuous. Therefore, due to Theorem 67 we get the
statement of this corollary.

\par {\bf 71. Definition.} A functional $J\in BC(G,{\sf g})^*$ is said
to have a compact support, if there exists a compact subset $V$ in
$G$ such that $J(f)=0$ for each $f\in BC(G,{\sf g})$ with $f(x)=0$
for every $x\in V$.
\par A hood on $G$ is a mapping $h: G\to [0,\infty )$ such that the set
$\{ x\in G: h(x)\ge \epsilon \} $ is compact for each $\epsilon
>0$. A subset $W\subset BC(G,{\sf g})$ is called strictly open if
for each $f\in W$ and a semi-norm $u\in \cal S$ in $\sf g$ there
exists a hood $h$ such that $W\supset \{ g\in BC(G,{\sf g}):
\sup_{x\in G} u(f(x)-g(x)) h(x) \le 1 \} $. Strictly open subsets in
$BC(G,{\sf g})$ form a topology in $BC(G,{\sf g})$ called the strict
topology.

\par {\bf 72. Theorem.} {\it The following conditions on $J\in
BC(G,Mat_n({\bf K}))^*$ are equivalent:
\par $(1)$ there exists $\mu \in M(G,Mat_n({\bf K}))$ such that $J=
\lambda _{\mu }$ (see Theorem 66);
\par $(2)$ for each $\epsilon >0$ there exists a compact subset
$V$ in $G$ such that $|J(f)| \le \max \{ \| J \| \sup_{x\in V} \|
f(x) \| , \epsilon \| f \| )$ for every $f\in BC(G,Mat_N({\bf K}))$;
\par $(3)$ $J$ is the limit of elements in $BC(G,Mat_n({\bf K}))^*$
having compact supports; \par $(4)$ $J$ is strictly continuous.}

\par {\bf Proof.} If $W_1$ and $W_2$ are strictly open, $f\in W=
W_1\cap W_2$, $u\in \cal S$, then there exist hoods $h_j$ such that
$W_j\supset \{ g\in BC(G,{\sf g}): \sup_{x\in G} u(f(x)-g(x)) h_j(x)
\le 1 \} $ for $j=1, 2$, then $W\supset \{ g\in BC(G,{\sf g}):
u(f(x)-g(x)) h(x) \le 1 \} $, where $h(x) = \max (h_1(x), h_2(x))$
for each $x$ is the hood such  that $ \{ x\in G: h(x)\ge \epsilon \}
= \{ x\in G: h_1(x)\ge \epsilon \} \cup \{ x\in G: h_2(x)\ge
\epsilon \}$ is compact for each $\epsilon >0$ as the union of two
compact sets. Thus strictly open subsets form a topology.

\par The $\bf K$-algebra $Mat_n({\bf K})$ is normed. Since $Mat_n({\bf
K})$ is the finite dimensional space over $\bf K$ its topologically
dual space is isomorphic with $Mat_n({\bf K})$. For a compact subset
$W$ in $G$ let $R_W$ denotes the restriction mapping $R_W:
BC(G,Mat_n({\bf K}))\to C(W,Mat_n({\bf K}))$. In view of Theorem
5.24 \cite{roo} there exists a $\bf K$-linear isometric embedding
$T_W: C(W,Mat_n({\bf K}))\hookrightarrow PC(G,Mat_n({\bf K}))$ such
that $R_W\circ T_W=I$, since $n\in \bf N$, where $PC(G,X)$ denotes
the closed $\bf K$-linear hull in $BC(G,X)$ of the subset $\{ Ch_A:
A\in {\sf Bco}(G), A \mbox{ is compact} \} $.
\par If $J=\lambda _{\mu }$ with $\mu \in M(G,Mat_n({\bf K}))$, then
$N_{\mu }$ is the hood and $|J(f)|\le \| f \|_{N_{\mu }}$ for every
$f\in BC(G,Mat_n({\bf K}))$, consequently, $J$ is strictly
continuous, that is, $(1)\Rightarrow (4)$.
\par If $J$ is strictly continuous take $\pi \in \bf K$ with $0<|\pi
| <1$. There exists a hood $h$ for which $ \{ f: \| f \|_h<1 \}
\subset \{ f: |J(f)|\le |\pi | \} $. Therefore, $|J(f)|\le \| f
\|_h$ for each $f$. For $\epsilon >0$ put $W := \{ x\in G: h(x)\ge
\epsilon \} $. If $f\in BC(G,Mat_n({\bf K}))$, then take $g :=
T_WR_Wf$, hence $J(f) = J(f-g) +J(g)$ and $|J(f-g)| \le \sup_{x\in
G} \| f(x)-g(x) \| h(x) \le \sup_{x\in G\setminus W} \| f(x) -g(x)
\| \epsilon \le \| f \| \epsilon $ and $|J(g)|\le \| J \| \| g \|
\le \| J \| \sup_{x\in W} \| f(x) \| $, consequently,
$(4)\Rightarrow (3)$.
\par Suppose that $(2)$ is satisfied, then $J_W$ has the compact
support. Therefore, $|J(f) - J_W(f)| = |J(f-T_WR_Wf)|\le \epsilon \|
f- T_WR_Wf \| <\epsilon \| f \| $, consequently, $ \| J -J_W \| \le
\epsilon $, hence $(2)\Rightarrow (3)$.
\par Let $(3)$ be satisfied. We have that $M(G,Mat_n({\bf K}))$ is
complete and $\lambda $ is the isometry, consequently, the range of
$\lambda $ is closed in $BC(G,Mat_n({\bf K}))^*$. Therefore, without
loss of generality consider $J$ with the compact support. Suppose
that $W\subset G$ and $J(f)$ for $f$ vanishing identically on $W$.
Put $\mu ^J_{i,j}(A) := J(E_{i,j}~ Ch_A)$ for each $A\in {\sf
Bco}(G)$ and all $i, j=1,...,n$, hence $\mu : {\sf Bco}(G)\to
Mat_n({\bf K})$ is additive, where $E_{i,j}$ is the matrix with the
element $1$ at the $(i,j)$-th place and zeros at others places, $\mu
(A)= \mu ^J(A)$ is the matrix with matrix elements $\mu
^J_{i,j}(A)$. Then $N_{\mu }(x)=0$ for every $x\in G\setminus W$ and
$N_{\mu }$ is bounded on $G$, consequently, $\mu $ is the measure on
${\sf Bco}(G)$.
\par The normed space $C(W,Mat_n({\bf K}))$ has an orthonormal base
consisting of functions $E_{i,j}Ch_A$, where $i, j=1,...,n$ and
$A\in {\sf Bco}(G)$. Suppose that $f\in BC(G,Mat_n({\bf K}))$, then
there exists $A_k\in {\sf Bco}(G)$ and $b_k\in Mat_n({\bf K})$ such
that $b_k = a_k E_{i(k),j(k)}$ and $\lim_{k\to \infty } a_k =0$ and
$f=\sum_k b_k Ch_{A_k}$ uniformly on $W$. Therefore, $J(f) =
J(\sum_k b_k Ch_{A_k}) = \sum_k a_k \mu _{i(k),j(k)} (A_k) =\int_G
\sum_k a_kCh_{A_k} \mu _{i(k),j(k)}(dx) = \int_G \sum_k b_k\mu (A_k)
= \int_G f(x)\mu (dx)$, since $\mu (A) = \sum_{i, j=1}^n E_{i,j}
J(E_{i,j}~ Ch_A)$, consequently, $(3)\Rightarrow (1)$.

\par {\bf 73. Corollary.} {\it Let $G$ and $H$ be zero-dimensional
Hausdorff spaces, let also $X$ be a complete topological algebra
over $\bf K$. If $\mu \in M(G,X)$ and $\nu \in M(H,X)$ are tight
measures, then $\mu \times \nu $ is a tight measure on $G\times H$,
$\mu \times \nu \in M(G\times H,X)$.}
\par {\bf Proof.} This follows from Theorem 72$(2)$.

\par {\bf 74. Example.}  Consider convolutions of tight measures.
Suppose that $G$ is a zero-dimensional Hausdorff topological
semi-group and $X$ be a topological algebra over $\bf K$. For $\mu
,\nu \in M(G,X)$ and $f\in BC(G,{\sf g})$ define $Jf := \int_G f(xy)
(\mu (dx) \times \nu (dy))$. If $u\in \cal S$ is a consistent
semi-norm in $X$ and $\sf g$, then $u(Jf) \le \sup_{x, y\in G}
u(f(xy)) N_{\mu ,u}(x) N_{\nu ,u}(y)$. For each $\epsilon >0$ we
have that $G_{\mu ,u,\epsilon } := \{ x\in G: N_{\mu ,u}(x) \ge
\epsilon \} $ and $G_{\nu ,u,\epsilon } := \{ y\in G: N_{\nu ,u
}(y)\ge \epsilon \} $ are compact, hence their product $G_{\mu
,u,\epsilon } G_{\nu ,u,\epsilon }$ is compact, moreover, $u(Jf) \le
\max \{ \sup \{ u(f(z)): z\in G_{\mu ,u,\epsilon } G_{\nu
,u,\epsilon } \} \| \mu \|_u \| \nu \|_u; \| f \|_u \| \mu \|_u
\epsilon ; \| f \|_u \| \nu \|_u\epsilon \} $. Therefore, $J$ is
induced by a tight measure denoted by $\mu *\nu $ such that $\int_G
f(x)[\mu *\nu ](dx) = \int_G f(xy) (\mu (dx) \times \nu (dy))$ for
each $f\in BC(G,X)$. In particular, for $f=Ch_A$ with $A\in {\sf
Bco}(G)$ we get \par $(i)$ $[\mu *\nu ](A) = (\mu \times \nu ) ( \{
(x,y)\in G\times G, xy\in A \} )$.\\ The tight measure $\mu *\nu $
is called the convolution product of $\mu $ and $\nu $. Evidently,
$(a\mu + b\zeta ) *\nu = (a\mu * \nu ) + (b\zeta * \nu )$ and $\mu *
(a\nu + b\zeta ) = (a\mu *\nu ) + (b\mu * \zeta )$ for each $a, b\in
\bf K$ and $\mu , \nu , \zeta \in M(G,X)$, since $\bf K$ is the
commutative field. Hence $M(G,X)$ is the algebra with the addition
$(\mu + \nu )(A) = \mu (A) +\nu (A)$ for each $A\in {\sf Bco}(G)$
and the multiplication given by the convolution product of measures.
\par From $(i)$ it follows, that $N_{\mu *\nu ,u}(z) = \sup_{x, y\in
G, xy=z} N_{\mu ,u}(x) N_{\nu ,u}(y)<\infty $, consequently,
$M(G,X)$ is the topological algebra with the family of semi-norms
$\| \mu \|_u := \sup_{x\in G} N_{\mu ,u}(x)$, $u\in \cal S$, such
that $\| \mu *\nu \|_u \le \| \mu \|_u \| \nu \|_u$ for each $\mu ,
\nu \in M(G,X)$.

\par {\bf 75. Lemma.} {\it The mapping 39$(SI)$
and Conditions 35$(M1-M4)$ induce an isometry between $L^2({\cal
R}(G),{\sf g})$ and $L^2(\xi ,{\sf g})$.}
\par {\bf Proof.} At first demonstrate, that there exists a
linear isometric mapping of $L^0({\cal R},{\sf g})$ on $L^0(\xi
,{\sf g})$. Let $f(x) =\sum_k a_k Ch_{A_k}(x)$ and $g(x)=\sum_l b_l
Ch_{A_l}(x)$ be simple functions in $L^0({\cal R},{\sf g})$, where
$a_k, b_l\in \sf g$ (see Definition 39). Then due to Conditions
35$(M1-M4)$ and 39$(1-7)$ there are the equalities: \par $(1)$ $M
[(\int_G [f(x)\xi (dx)), (\int_G g(x) \xi (dx))] = \sum_k [a_k,b_k]
\mu (A_k)$ \par $ =\int_G [f(x),g(x)] \mu (dx)$, \\
since ${\sf g}^2\ni \{ a,b \} \mapsto [a,b]\in Lc({\sf g})$ is the
continuous mapping, $\mu (A)\in Lc(X)$ for each $A\in {\cal R}$.
\par In view of Lemma 38 there exists the $\bf K$-valued
measure $Tr \mu $. Therefore,
\par $(2)$ $M ((\int_G f(x)\xi (dx)), (\int_G g(x) \xi (dx)))
= \sum_k Tr (b_k^Ta_k\mu )(A_k)$
\par $ =\int_G Tr (g^T(x)f(x) \mu )(dx)$, \\
since ${\sf g}^2\ni \{ a,b \} \mapsto (a,b)\in \bf K$ is the
continuous mapping from ${\sf g}^2$ into the field $\bf K$, where
$(Tr \mu )(A) := Tr \mu (A)$ for each $A\in {\cal R}(G)$. \par If
$F_1\in Lin (X)$ and $F_2\in Lc (X)$, then $F_1F_2\in Lc (X)$. To
each $b\in \sf g$ there corresponds the $\bf K$-linear continuous
operator $X\ni x\mapsto bx\in X$. Naturally $Lin (X)$ and $Lc (X)$
are the left $\sf g$-modules, since $(bF)\in Lin (X)$ for each $F\in
Lin (X)$ and $(bF)\in Lc(X)$ for every $F\in Lc (X)$ and each $b\in
\sf g$, where $(bF)(x) := b(Fx)$ for all $x\in X$. Since $\mu (A)\in
Lc(X)$ and $g^T(x)f(x)\in \sf g$, then $g^T(x)f(x)\mu (A)\in Lc (X)$
for each $x\in G$ and $A\in {\cal R}(G)$ and all $f, g\in L^0({\cal
R},{\sf g})$.
\par Take in $L^0({\cal R},{\sf g})$ the semi-norm
\par $(3)$ $\| f\|_{2,  \mu } =
[\sup_{x\in G} N_{Tr f^Tf\mu  }(x)]$ \\
and in $L^0(\xi ,{\sf g})$ put
\par $(4)$ $\| \eta _f\|_{2,N_P} :=
[\sup_{x\in G} |(f(x)\xi (x), f(x)\xi (x))|N_P(x)]^{1/2}$ \\
for each $\eta = \int_G f(x)\xi (dx)$. The semi-norm $(4)$ is
continuous relative to the family of semi-norms 43$(1)$.
\par Semi-norm $(3)$ is continuous relative to the family of semi-norms:
\par $(5)$ $\| f \|_{\mu ,u} := \sup_{x\in G} N_{f^Tf\mu ,u}(x)$.
\par Since $M((a\xi
(A),b\xi (A)))=Tr (b^Ta\mu (A))$ for each $A\in {\cal R}(G)$, then
$N_{Tr f^Tf\mu }(x) = \inf_{A\in {\cal R}(G), x\in A} \| A \|_{Tr
f^Tf\mu }$, where $(f^Tf\mu )(dx) = f^T(x)f(x)\mu (dx)$. At the same
time $M(a\xi (B),b\xi (B)) = \int_{\Omega }(a\xi (\omega ,B), b\xi
(\omega ,B)) P(d\omega )$, $|M(a\xi (B),b\xi (B))|\le \sup_{\omega
\in \Omega } |(a\xi (\omega ,B),b\xi (\omega ,B))| N_P(\omega )$ for
each $a, b\in \sf g$.
\par In view of Lemma 38 $Tr g^Tf\mu $ is the measure for each $f, g\in
L^0({\cal R},{\sf g})$, consequently, taking a shrinking family
$\cal S$ in ${\cal R}(G)$ such that $\bigcap_{A\in \cal S} A= \{ x
\} $ gives
\par $N_{Tr g^Tf\mu }(x)=
\inf_{A\in {\cal R}(G), x\in A}[\sup_{B\in {\cal R }(G), B\subset
A_k, k } \sup_{\omega \in \Omega , k} |(a_k\xi (\omega ,B), b_k\xi
(\omega ,B))| N_P(\omega )]$. \\ Thus $N_{Tr g^Tf\mu }(x)=\inf_{A\in
{\cal R}(G), x\in A} [\sup_{B\in {\cal R}(G), B\subset A_k, k } \|
(a_k\xi (*,B), b_k\xi (*,B)) \|_{L^2(P,{\bf K})}]$ and $\| A_k
\|_{Tr g^Tf\mu } = \| (a_k\xi (*,A_k), b_k\xi (*,A_k)) \|_{L^2(P)}]$
for each $k=1,...,m$ due to Lemma 2 and due to the choice $ \| A_k
\| _{Tr b_k^Ta_k\mu } = |Tr b_k^Ta_k\mu (A_k) |$ without loss of
generality.
\par On the other hand, $M[a\xi (A),b\xi (A)]=b^Ta\mu (A)\in Lc(X)$
for each $A\in {\cal R}(G)$. Then $N_{b^Ta\mu ,u}(x) = \inf_{A\in
{\cal R}(G), x\in A}[\sup_{B\in {\cal R }(G), B\subset A} u(b^Ta\mu
(B))]$ \\  $\le u(b^Ta)N_{\mu ,u}(x)<\infty $, where $N_{\mu ,u}(x)
:= \inf_{A\in {\cal R}(G), x\in A}[\sup_{B\in {\cal R }(G), B\subset
A} u(\mu (B))$. Therefore, $\| f \|_{\mu ,u}= \| f \|_{2,P,u}$ for
each $f\in L^0({\cal R},{\sf g})$ and each consistent semi-norm $u$
in $\sf g$, $X$ and $Lin (X)$.
\par The mapping $\psi $ from 39$(SI)$ also is $\bf K$-linear
from $L^0(\xi ,{\sf g})$ into $L^0({\cal R}(G),{\sf g})$ such that
$\psi $ is the isometry relative to the consistent semi-norms
43$(1)$ and 75$(5)$ due to Formula $(1)$ and Lemmas 2, 45 and
Theorem 56.
\par Two spaces $L^2(P,{\sf g})$ and $L^2(\mu ,{\sf g})$ are
complete by their definitions, consequently, $\psi $ has the $\bf
K$-linear extension from $L^2({\cal R}(G),{\sf g})$ onto $L^2(\xi
,{\sf g})$ which is the isometry between $L^2({\cal R}(G),{\sf g})$
and $L^2(\xi ,{\sf g})$.

\par {\bf 76. Definition.} If $f\in L^2({\cal R}(G),{\sf g})$,
then put by the definition:
\par $\eta =\psi (f)=\int_G f(x)\xi (dx)$.
\par The random vector $\eta $ we call the non-archimedean
stochastic integral of the function $f$ by an orthogonal stochastic
measure $\xi $.

\par {\bf 77. Remark.} Consider random vectors of the form:
\par $\eta (t) = \int_G g(t,x)\xi (dx)$, where $\xi $ is an
orthogonal stochastic measure on a measurable space $(G,{\cal
R}(G))$ with values in $X$ and a structural measure $\mu $ with
values in $Lc(X)$ as above, $t\in T$, $g(t,x)\in L^2(G,{\cal
R}(G),\mu ,{\sf g})$ as the function by $x\in G$ for each $t\in T$,
where $T$ is a set.

\par The covariance operator of a random vector $\eta $ is
\par $(1)$ $B(t_1,t_2) = M\{ \eta ^T(t_1), \eta (t_2) \} = \int_G
\{ g^T(t_1,x), g(t_2,x) \} \mu (dx)$, moreover,
\par $(2)$ $M\{ \eta (t_1), \eta ^T(t_2) \} =
\int_G \{ g(t_1,x), g^T(t_2,x) \} Tr \mu (dx)$ \\
in the notation of \S \S 31, 39, where $X$ is the left $\sf
g$-module, while $Lc(X)$ is supplied with the natural structure of
the left $Lc({\sf g})$-module, $ \{ a^T, b \} \in Lin ({\sf g})$,
$\{ a, b^T \} \in \sf g$ for each $a, b\in \sf g$. Denote by $L^2 \{
g \} $ the closure in $L^2(G,{\cal R}(G),\mu ,{\sf g})$ of the $\bf
K$-linear span of the family of functions $ \{ g(t,x): t\in T \} $.
Therefore, $L^2 \{ g \} $ is the $\bf K$-linear closed subspace in
$L^2(G,{\cal R}(G),\mu ,{\sf g})$. If $L^2 \{ g \} = L^2(G,{\cal
R}(G),\mu ,{\sf g})$, then the system of functions $ \{ g(t,x): t\in
T \} $ is called complete in $L^2(G,{\cal R}(G),\mu ,{\sf g})$.

\par Let $\{ \eta (t): t\in T \} $ be a $X$-valued random vector.
Denote by $L^0 \{ \eta \} $ the family of all random vectors of the
form:
\par $\zeta = \sum_{k=1}^l a_k\eta (t_k)$, where $l\in \bf N$,
$t_k\in T$, $a_k\in \sf g$. Then $L^2 \{ \eta \} $ denotes the
closure of $L^0 \{ \eta \} $ in $L^2(\Omega ,{\cal R},P,X)$.
\par A family of random vectors $ \{ \zeta _{\beta }: \zeta _{\beta }
\in L^2(\Omega ,{\cal R},P,X); \beta \in \Lambda \} $ is called
subordinated to the random $X$-valued function $\{ \eta (t): t\in T
\} $, if $\zeta _{\beta }\in L^2 \{ \eta \} $ for each $\beta \in
\Lambda $.

\par {\bf 78. Lemma.} {\it If $X$ is a left $\sf g$-module, then
$Lin(X)$ is a left $Lin ({\sf g})$-module, $Lc(X)$ is a left $Lin
({\sf g})$-module as well as left $Lc({\sf g})$-module.}
\par {\bf Proof.} Let $x\in X$, $y\in \sf g$, $A\in Lin ({\sf g})$
and $B\in Lin (X)$, then $Bx\in X$, $y(Bx)\in X$ and $Ay \in \sf g$,
hence $AyBx\in X$, since the left module is associative.
Particularly for $1\in \sf g$ it gives $AB\in Lin (X)$. Therefore,
there exists the multiplication $Lin ({\sf g})\times Lin (X)\ni
(A,B)\mapsto AB\in Lin(X)$. If $B\in Lc(X)$, then $AyB\in Lc(X)$,
since $vB\in Lc(X)$ for each $v\in \sf g$. In particular, for $y=1$,
consequently, $Lc(X)$ is the left $Lin({\sf g})$-module and
inevitably left $Lc({\sf g})$-module, since $Lc({\sf g})\subset
Lin({\sf g})$.

\par {\bf 79. Theorem.} {\it Let a covariance operator $B(t_1,t_2)$
of a random $X$-valued function $ \{ \eta (t): t\in T \} $ admits
representation 76$(1)$, where $X$ and $\sf g$ are as in \S \S 30 and
39, $\mu $ is a $Lc(X)$-valued measure on $(G, {\cal R}(G))$, $\mu
^T=\mu $, $g(t,x)\in L^2(G,{\cal R}(G),\mu ,Lc(X);{\sf g})$ for each
$t\in T$ and the family $ \{ g(t,x): t\in T; \} $ is complete in
$L^2(G,{\cal R}(G),\mu ,Lc(X);{\sf g})$. Then $\eta (t)$ can be
presented
in the form: \par $(1)$ $\eta (t) = \int_G g(t,x)\xi (dx)$ \\
with probability $1$ for each $t\in T$, where $\xi $ is a stochastic
orthogonal $X$-valued measure subordinated to the random function
$\eta (t)$ and with a structure function $\mu $.}

\par {\bf Proof.} Consider functions of the form:
\par $(2)$ $f(x) = \sum_{k=1}^l b_k g(t_k,x)$, \\
where $t_k\in T$, $b_k\in \sf g$, $l\in \bf N$. Put
\par $(3)$ $\psi (f) =\zeta = \sum_{k=1}^l b_k\eta (t_k)$.
Denote by $L^0 \{ g \} $ the family of all vectors of Form $(2)$. In
$L^0 \{ g \} $ there is the $\bf K$-bi-linear functional:
\par $(4)$ $(f_1,f_2) := \int_G \{ f_1(x), f_2^T(x) \} Tr \mu (dx)$.
\\
In view of Lemma 75 the mapping $\zeta = \psi (f)$ is the $\bf
K$-linear topological isomorphism of $L^0 \{ g \} $ onto $L^0 \{
\eta \} $. When particularly $X$ and $\sf g$ are normed, then $\psi
$ is the isometry. Thus $\psi $ has the continuous extension up to
the $\bf K$-linear topological isomorphism of $L^2 \{ g \} $ onto
$L^2 \{ \eta \} $. \par If $A\in {\cal R}(G)$, then $Ch_A\in
L^2(G,{\cal R}(G),\mu ,{\sf g})$, since $1\in \sf g$. But $L^2 \{ g
\} = L^2(G,{\cal R}(G),\mu ,{\sf g})$ due to completeness of the
family $ \{ g(t,x): t\in T \} $. Therefore, $Ch_A\in L^2(G,{\cal
R}(G),\mu ,{\sf g})$. Put $\xi (A) := \psi (Ch_A)$, then $\xi (A)$
is the orthogonal stochastic measure with the structure function
$\mu $ due to Lemma 78, since $(5)$ $M\{ \xi ^T(A), \xi (B) \} =
\int_G \{ Ch_A^T(x), Ch_B(x) \} \mu (dx)= \mu (A\cap B)$  for each
$A, B\in {\cal R}(G)$.
\par Let now $\gamma (t) := \int_G g(t,x) \xi (dx)$.
Since $M \{ \eta ^T(t), \xi (A) \} = \int_G \{ g^T(t,x), Ch_A(x) \}
\mu (dx)$ and $M \{ \xi ^T(A), \eta (t) \} = \int_G \{ Ch_A^T(x),
g(t,x) \} \mu (dx)$ and $\psi $ is the $\bf K$-linear topological
isomorphism, then $M \{ \eta ^T(t), \gamma (t) \} =M \{ \gamma
^T(t), \eta (t) \} = \int_G \{ g^T(t,x), g(t,x) \} \mu (dx)$.
Therefore, $M \{ Ch_A(\eta (t) - \gamma (t))^T, Ch_A(\eta (t) -
\gamma (t)) \} = M \{ Ch_A\eta ^T(t), Ch_A\eta (t) \} - M \{
Ch_A\eta ^T, Ch_A\gamma (t) \} - M \{ Ch_A\gamma ^T(t), Ch_A\eta (t)
\} + M \{ Ch_A\gamma ^T(t), Ch_A\gamma (t) \} =0$ for each $A\in
{\cal R}(G)$, consequently, 79$(1)$ is accomplished with probability
$1$ for each $t\in T$.

\par {\bf 80. Definition.} Let $\eta (t)$ be a $\sf g$-valued stochastic
process or stochastic function such that for each $n\in \bf N$ and
each $t_1,...,t_n\in \bf T$ with $t, t+t_1, ..., t+t_n\in \bf T$ the
mutual distribution of $\eta (t+t_1), ..., \eta (t+t_n)$ is
independent from $t$, where $\bf T$ is an additive semigroup. Then
$\eta (t)$ is called the stationary stochastic function, where $P:
{\cal R}\to \bf K$ is a probability measure.

\par Suppose that $\bf T$ is a uniform space.
A $\sf g$-valued stochastic function $\eta (t)\in L^b(\Omega ,{\cal
R},P,{\bf K};{\sf g})$, $t\in \bf T$, $1\le b<\infty $, is called
mean-$b$-continuous at $t_0\in \bf T$, if there exists $\lim_{t\to
t_0} \eta (t) =\eta (t_0)$ in the sense of convergence in the space
$L^b(\Omega ,{\cal R},P,{\bf K};{\sf g})$, when $t$ tends to $t_0$
in $\bf T$. In particular, for $b=2$ it is mean-square continuity
and convergence respectively. If $\eta (t)$ is mean-$b$-continuous
at each point of $\bf T$, then $\eta (t)$ is called
mean-$b$-continuous on $\bf T$.

\par {\bf 81.} Suppose that an algebra $\sf g$ over $\bf C_p$ has
a uniformity $\tau $ relative to which it is complete. Let $\sf g$
has a $\bf Q_p$-linear embedding into $c_0(\gamma ,{\bf Q_p})$ for
some set $\gamma $ such that the norm uniformity $n_u$ in $\sf g$
inherited from the Banach space $c_0(\gamma ,{\bf Q_p})$ with the
standard norm $ \| * \| $ is not stronger, than $\tau $, that is
$n_u\subset \tau $, moreover, $\sf g$ is everywhere dense in
$(c_0(\gamma ,{\bf Q_p}), \| * \| )$.

\par {\bf Theorem.} {\it Let $\eta (t)$ be a stationary
mean-square-continuous stochastic process with values in $\sf g$,
$t\in {\bf T} = \bf C_r$, where $r$ and $p$ are mutually prime
numbers, $M\eta (t)=0$, then there exists an orthogonal $\sf
g$-valued stochastic measure $\xi (A)$ on ${\sf Bco}({\bf C_r})$
subordinated to $\eta (t)$ such that
\par $(1)$ $\eta (t) = \int_{{\bf C_r}} g(t,x) \xi (dx)$, \\
where $g(t,x)$ is a $\bf C_p$-valued character from the additive
group $({\bf C_r}, +)$ into the multiplicative group $({\bf
C_p},\times )$. Between $L^2 \{ \eta \} $ and $L^2 \{ \mu \} $ there
exists a $\bf K$-linear topological isomorphism $\psi $ such that
\par $(2)$ $\psi (\eta (t)) = g(t,*)$, $\psi (\xi (A)) = Ch_A$, if
\par $(3)$ $\zeta _j =\psi (f_j)$, then $\zeta _j = \int_{\bf C_r}
f_j(x)\xi (dx)$ and $M \{ \zeta _1^T, \zeta _2 \} = \int_{\bf C_r}
\{ f_1(x)^T, f_2(x) \} \mu (dx)$.}

\par {\bf 82. Definition.} Formula 81$(1)$ is called the spectral
decomposition of the stationary stochastic process. A measure $\xi
(A)$ is called a stochastic spectral measure of the stationary
stochastic process $\eta (t)$.

\par {\bf Proof of Theorem 81.} Since $\eta (t)$ is a stationary stochastic
process, then for each continuous function $f: {\sf g}^n\to \bf K$
the mean value $Mf(\eta (t+t_1),...,\eta (t+t_n))$ is independent
from $t$, where $n\in \bf N$. By the condition of this theorem $\eta
(t)\in L^2(\Omega ,{\cal R},P,{\bf C_p};{\sf g})$, consequently,
there exist $M\eta (t)=m$ and
\par $(4)$ $M\{ [\xi (t)-m]^T, [\xi (q)-m] \} =B(t-q)\in Lc ({\sf g})$
for each $t, q\in \bf C_r$, where $m=0$. Evidently $B^T(t-q) =
B(q-t)$ for each $t, q\in \bf C_r$.

\par Consider now ${\bf C_p}$-valued characters of $({\bf C_r}, +)$
as the additive group, where $r=p'$, $p\ne p'$ are prime numbers.
For $p$-adic numbers $x = \sum_{k=N}^{\infty }x_kp^k$, where $x\in
\bf Q_p$, $x_k\in \{ 0, 1,...,p-1 \} $, $N\in \bf Z$, $N=N(x)$,
$x_N\ne 0$, $x_j=0$ for each $j<N$, put as usually $ord_p (x) = N$
for the order of $x$, thus its norm is $|x|_{\bf Q_p} =p^{-N}$.
Define the function $[x]_{\bf Q_p} :=\sum_{k=N}^{-1}x_kp^k$ for
$N<0$, $[x]_{\bf Q_p}=0$ for $N\ge 0$ on ${\bf Q_p}$. Therefore, the
function $[x]_{\bf Q_p}$ on $\bf Q_p$ is considered with values in
the segment $[0,1]\subset \bf R$.
\par Consider the field $\bf C_r$ as the vector space over the field
$\bf Q_r$. There is a multiplicative non-archimedean norm $|
*|_{\bf C_r} = |*|$ in $\bf C_r$, which gives the uniformity in it.
Take an equivalent uniformity given by a norm $|*|_r$ such that
$|x|_r\in \{ r^l: l\in {\bf Z} \} \cup \{ 0 \} $ for each $x\in \bf
C_r$. If $x\ne 0$ put $|x|_r := \min \{ r^l: |x|_{\bf C_r}\le r^l,
l\in {\bf Z} \} $, $|0|_r = 0$, hence \par $(i)$ $|x|_r/r \le
|x|_{\bf C_r} \le |x|_r$ for each $x\in \bf C_r$ \\ and inevitably
$\bf C_r$ is the topological vector space over $\bf Q_r$ relative to
$|x|_r$. Since $\bf C_r$ is the extension of $\bf Q_r$, then the
restriction of $|*|_r$ on $\bf Q_r$ is the $r$-adic norm. On the
entire $\bf C_r$ this $|*|_r$ in general need not be multiplicative.
Verify, that it is indeed a non-archimedean norm. At first $|x|_r\ge
0$ for each $x\in \bf C_r$, $|x|_r=0$ if and only if $x=0$ due to
$(i)$. If $x, y\in \bf C_r$, $x\ne 0$, $y\ne 0$, then $|x|=r^a$,
$|y| =r^b$, $|x+y|=r^c$ with $a, b, c\in \bf R$, $c\le \max (a, b)$,
where $r\ge 2$ is a prime number. Then $|x|_r = r^A$, $|y|_r = r^B$,
$|x+y|_r = r^C$, where $a\le A$, $b\le B$, $c\le C$, $A, B, C\in \bf
Z$ are the least integers satisfying these inequalities. Therefore,
$C\le max (A,B)$, consequently, $|x+y|_r\le \max (|x|_r, |y|_r)$ for
each $x, y\in \bf C_r$. \par In view of Theorems 5.13 and 5.16
\cite{roo} the $\bf Q_r$-linear space $({\bf C_r}, |*|_r)$ is
isomorphic with $c_0(\alpha ,{\bf Q_r})$, where $\alpha $ is a set,
which is convenient to consider as an ordinal due to Zermelo theorem
\cite{eng}.
\par Let $(x,y) := (x,y)_{\bf Q_r} := \sum_{j\in \alpha }x_jy_j$ for
$x, y\in \bf C_r$, $x = (x_j: j\in \alpha , x_j\in {\bf Q_r} )$.
This series $(x,y)$ converges in $\bf Q_r$, since for each $\epsilon
>0$ the set $ \{ j: |x_j|_r\ge \epsilon \} $ is finite.

\par If $X$ is a complete locally $\bf C_r$-convex space, then it is
the projective limit of Banach spaces $V_u := X/Y_u$ over $\bf C_r$,
where $Y_u := \{ x\in X: u(x)=0 \} $, $u$ is a semi-norm in $X$,
$u\in \cal S$ \cite{nari}. Each $V_u$ can be supplied with the
structure of a Banach space over $\bf Q_r$. Therefore, $X$ can be
supplied with the structure $X_r$ of the complete locally $\bf
Q_r$-convex space with a topology $\tau _r$.

\par Consider the case of such $X$, when $X_r$ has an embedding into
$c_0 (\beta , {\bf Q_r})$ for some $\beta \ge \alpha $ and the norm
topology $n_r$ of $|*|_r$ in $X_r$ inherited from $c_0(\beta ,{\bf
Q_r})$ is such that $\tau _r\supset n_r$. Then each $\bf Q_r$-linear
continuous functional on $(c_0(\beta ,{\bf C_r}), |*|_r)$ is also
continuous on $(X_r, \tau _r)$.

\par Define also a character with values in $\bf C_p$ for $(X,+)$
as the additive group, $r\ne p$. Put
\par $\chi _{r,p; s} (x)=\epsilon ^{ [(s,z)_{\bf Q_r}]_{\bf Q_r}/z}$, where
$\epsilon =1^z$ is a root of unity in $\bf C_p$, $z=r^{ord_r
[(s,z)_{\bf Q_r}]_{\bf Q_r}}$, $s, z\in X_r$ or we can consider $s,
z$ as elements in $X$ as well (see above).

\par For a tight measure $\mu : {\cal R}(X)\to Lc ({\sf g})$
or $\mu : {\cal R}(X)\to \bf C_p$ the characteristic functional
$\hat \mu $ is given by the formula: ${\hat \mu }(s) := \int_X \chi
_{r,p;s} (z)\mu (dz)$, where $s\in X_r$, $X$ is over $\bf C_r$.
\par In general the characteristic functional of the measure
$\mu : {\cal R}(G)\to Lc({\sf g})$ or $\mu : {\cal R}(G)\to \bf C_p$
is defined in the space $C^0(G,{\bf C_r})$ of continuous functions
$f: G\to \bf C_r$
\par ${\hat \mu }(f) := \int_G \chi _{r,p;1} (f(z))\mu (dz)$,
where $1\in \bf C_r$, $G$ is a totally disconnected topological
Hausdorff space with a covering ring ${\cal R}(G)$.

\par In view of Theorems 2.21 and 2.30 \cite{lujmsqim} and Theorem 64
and $(4)$ above there exists a $Lc ({\sf g})$-valued measure $\mu $
on ${\sf Bco}({\bf C_r})$ such that $B(t) = \int_{\bf C_r} \chi
_{r,p; 1} (ty) \mu (dy)$. Functions $g(t,y):= \chi _{r,p; 1} (ty)$
are continuous and uniformly bounded. Since $|g(t,y)|_{C_p}=1$ for
each $t, y\in \bf C_r$, then $g(t,y)\in L^2 \{ \mu \} $. \par In
view of the Kaplansky Theorem A.4 \cite{sch1} and Theorem 56 above
the family of functions $ \{ g(t,y): t, y \in {\bf C_r} \} $ is
complete in $L^2({\sf g},{\sf Bco}({\sf g}),\mu ,Lc({\sf g});{\bf
C_p})$.  Thus statements 81$(1-3)$ follow from Theorem 79.

\end{document}